\documentclass[12pt]{amsart}

\usepackage{amsmath}
\usepackage{amscd}
\usepackage{pb-diagram}
\usepackage{comment}
\usepackage{graphicx}
\usepackage{amssymb}
\usepackage{lmodern}
\usepackage{color}
\usepackage{wrapfig}
\let\visiblecomments y
\def\refeq#1{\if\workingver y(\ref{#1})-[[#1]]\else(\ref{#1})\fi}
\def\refth#1{\if\workingver y\ref{#1}-[[#1]]\else\ref{#1}\fi}
\def\mylabel#1{\if\workingver y\label{#1}{\bf\ \ [[#1]]\ \ }\else\label{#1}\fi}
\def\mybibitem#1{\if\workingver y\bibitem{#1}{\bf\ \ [[#1]]\ \
}\else\bibitem{#1}\fi}

\newfont{\msam}{msam10}
\newfont{\msbm}{msbm10}

\def\cA{\text{$\mathcal A$}}
\def\cB{\text{$\mathcal B$}}

\def\cM{\text{$\mathcal M$}}

\def\cO{\text{$\mathcal O$}}

\def\cS{\text{$\mathcal S$}}
\def\cT{\text{$\mathcal T$}}

\def\cV{\text{$\mathcal V$}}
\def\cW{\text{$\mathcal W$}}

\def\bA{\text{$\mathbf A$}}
\def\bB{\text{$\mathbf B$}}
\def\bC{\text{$\mathbf C$}}

\def\bP{\text{$\mathbf P$}}
\def\bQ{\text{$\mathbf Q$}}
\def\bR{\text{$\mathbf R$}}
\def\bS{\text{$\mathbf S$}}

\newcommand{\St}{{\operatorname{St}\,}}

\newcommand{\id}{\operatorname{id}}
\newcommand{\cl}{\operatorname{cl}}

\newcommand{\inte}{\operatorname{int}}

\newcommand{\dom}{\operatorname{dom}}

\renewcommand{\emptyset}{\varnothing}
\renewcommand{\subset}{\subseteq}

\newcommand{\bcap}{\mathbin{\bar{\cap}}}

\newcommand{\co}{\operatorname{co}}

\def\proof{{\bf Proof:\ }}

\def\begeq#1{\begin{equation}\mylabel{#1}}
\def\endeq{\end{equation}}

\def\mathobj#1{\mbox{$#1$}}

\def\NN{\mathobj{\mathbb{N}}}

\def\RR{\mathobj{\mathbb{R}}}

\def\ZZ{\mathobj{\mathbb{Z}}}

\def\setof#1{\mbox{$\{\,#1\,\}$}}

\def\cell#1{\protect\mbox{$\stackrel{\circ}{#1}$}}
\newcommand{\mto}{\multimap}
\newcommand{\pto}{\nrightarrow}
\newcommand{\union}[1]{\langle#1\rangle}

\newcommand{\famofgr}[1]{#1^{\shortparallel}}

\usepackage[all]{xy}
\usepackage{algorithm}
\usepackage[noend]{algpseudocode}

\newtheorem{theorem}{Theorem}[section]

\newtheorem{thm}[theorem]{Theorem}

\newtheorem{defn}[theorem]{Definition}
\newtheorem{cor}[theorem]{Corollary}
\newtheorem{prop}[theorem]{Proposition}

\definecolor{darkgreen}{rgb}{0.0, 0.8, 0.0}
\definecolor{darkred}{rgb}{1, 0.1, 0.3}
\definecolor{darkblue}{rgb}{0.1, 0.1, 1}

\graphicspath{{example-figs-100/}}

\title[Morse Decompositions in Combinatorial Dynamics]{Persistent Homology of Morse Decompositions in Combinatorial Dynamics}

\author[T.K.\ Dey]{Tamal K. Dey}
\address{Tamal K. Dey, Department of Computer Science and Engineering, The Ohio State University, Columbus, Ohio, USA}
\email{tamaldey@cse.ohio-state.edu}

\author[M.\ Juda]{Mateusz Juda}
\address{Mateusz Juda, Division of Computational Mathematics,
  Institute of Computer Science and Computational Mathematics,
  Faculty of Mathematics and Computer Science,
  Jagiellonian University, ul.~St. \L{}ojasiewicza 6, 30-348~Kra\-k\'ow, Poland.
}
\email{mateusz.juda@uj.edu.pl}

\author[T.\ Kapela]{Tomasz Kapela}
\address{Tomasz Kapela, Division of Computational Mathematics,
  Institute of Computer Science and Computational Mathematics,
  Faculty of Mathematics and Computer Science,
  Jagiellonian University, ul.~St. \L{}ojasiewicza 6, 30-348~Kra\-k\'ow, Poland.
}
\email{tomasz.kapela@uj.edu.pl}

\author[J.\ Kubica]{Jacek Kubica}
\address{Jacek Kubica, Division of Computational Mathematics,
  Institute of Computer Science and Computational Mathematics,
  Faculty of Mathematics and Computer Science,
  Jagiellonian University, ul.~St. \L{}ojasiewicza 6, 30-348~Kra\-k\'ow, Poland.
}
\email{jacek.kubica@uj.edu.pl}

\author[M.\ Lipi\'nski]{Micha\l{} Lipi\'nski}
\address{Micha\l{} Lipi\'nski, Division of Computational Mathematics,
  Institute of Computer Science and Computational Mathematics,
  Faculty of Mathematics and Computer Science,
  Jagiellonian University, ul.~St. \L{}ojasiewicza 6, 30-348~Kra\-k\'ow, Poland.
}
\email{michal.lipinski@uj.edu.pl}

\author[M.\ Mrozek]{Marian Mrozek}
\address{Marian Mrozek, Division of Computational Mathematics,
  Institute of Computer Science and Computational Mathematics,
  Faculty of Mathematics and Computer Science,
  Jagiellonian University, ul.~St. \L{}ojasiewicza 6, 30-348~Kra\-k\'ow, Poland.
}
\email{marian.mrozek@uj.edu.pl}
\thanks{TD is partially supported by NSF, USA Grants CCF 1740761 and CCF 1526513;
MJ, TK and MM are partially supported
       by the Polish National Science Center under Ma\-estro Grant 2014/14/A/ST1/00453.}

\date{}

\begin{document}

\begin{abstract}
  We investigate combinatorial dynamical systems on simplicial complexes considered as {\em finite topological spaces}.
  Such systems arise in a natural way from sampling dynamics and may be used to reconstruct some features
  of the dynamics directly from the sample.
  We study the homological persistence of {\em Morse decompositions} of such systems as a tool for validating the reconstruction.
  Our approach may be viewed as a step toward applying the classical persistence theory to data collected from a dynamical system.
  We present experimental results on two numerical examples.
\end{abstract}

\maketitle

\section{Introduction}
\label{sec:intro}

The aim of this research is to provide a tool for studying the topology of {\em Morse decompositions} of
sampled dynamics, that is dynamics known only from a sample.
Morse decomposition of the phase space of a dynamical system consists of a finite collection
of isolated invariant sets, called {\em Morse sets}, such that the dynamics outside the Morse sets
is gradient-like. This fundamental concept introduced in 1978 by Conley \cite{conley:78a}
generalizes classical Morse theory to non-gradient dynamics.
It has become an important tool in the study of the asymptotic behavior of flows, semi-flows
and multivalued flows (see \cite{BCCA2013,CoVa17,Souza2012} and the references therein).
Morse decompositions and the associated Conley-Morse graphs \cite{Ar09,BuAtAl2012}
provide a global descriptor of dynamics. This makes them an excellent object for studying the dynamics
of concrete systems. In particular, they have been recently applied in such areas
as mathematical epidemiology \cite{KnPiRo2015}, mathematical ecology \cite{Ar09,BuAtAl2012,GuLaSu2017} or
visualization \cite{Sz2013,WWL2016}.

Unlike the case of theoretical studies, the methods of classical mathematics do not suffice in most problems concerning concrete dynamics.
This is either because there is no analytic solution to the differential equation describing the system or, even worse,
the respective
equation is only vaguely known or not known at all.
In the first case the dynamics is usually studied by numerical experiments.
In some cases this may suffice to make mathematically rigorous claims about the system \cite{MiMr1995}.
In the latter case one can still get some insight into the dynamics by collecting data from physical experiments or observations,
for instance as a time series
\cite{ABMS2015,GuLaSu2017,MiMrReSz1999}.
In both cases the study is based on a finite and not precise sample, typically in the form of a data set.
The inaccuracy in the data may be caused by noise, experimental error, or numerical error.
Consequently, it may distort the information gathered from the data, raising the question whether the information is
trustworthy. One of possible remedies is to study the stability of the information with respect to perturbation of the data.
This approach to Morse decompositions constructed from samples is investigated in \cite{Sz2013} in the setting of
piecewise constant vector fields on triangulated manifold surfaces. The outcome of the algorithm proposed in \cite{Sz2013}
is the Morse merge tree which encodes the zero-dimensional persistence under perturbations of individual Morse sets
in the Morse decomposition.

In this paper we study general persistence of Morse decompositions in combinatorial dynamics,
not necessarily related to perturbations.
To this end, we define a homological persistence of the Morse decomposition over a sequence of combinatorial dynamics.
By combinatorial dynamics we mean a multivalued map acting on a simplicial complex treated as a finite topological space.
This general setting may be applied either to a finite sample of the action of a map on a subspace of $\RR^d$
\cite{BauerEdelsbrunnerJablonskiMrozek-2015,EdelsbrunnerJablonskiMrozek-2015} or to a combinatorial vector field \cite{Fo98} and its generalization multivector field \cite{Mr2016}.
The persistence is obtained by linking the homology of topologies induced by Morse decompositions and Alexandrov topology in the sequence of
combinatorial dynamical systems connected with continuous maps in zig-zag order.

On the theoretical level, the results presented in this paper may be generalized to arbitrary finite $T_0$ topological spaces.
From the viewpoint of applications, the finite topological
space may be
a collection of cells of a simplicial, cubical, or general cellular complex approximating a cloud of sampled points.
The multivalued map may be constructed either from the action of a given map on the set
of a sample points or from the available vectors of a sampled vector field.
The framework for persistence of Morse decompositions in the combinatorial setting developed in this
paper is general and may be applied to many different problems.

The language of finite topological spaces (see Section~\ref{ssec:fts})
enables us to emphasize differences between the classical and combinatorial dynamics.
These differences matter when the available data set is sparse and is difficult to be enriched.
In particular, in the classical setting the phase space has Hausdorff topology ($T_2$ topology) and
the Morse sets are compact. Hence, Morse sets are isolated since they are always disjoint.
To achieve such isolation in sampled dynamics, one needs data
not only in the Morse sets but also between the Morse sets. This may be a problem if
the available data set is sparse and cannot be enhanced.
Fortunately, the finite topological spaces in general are not $T_2$. Every set is compact but compactness
does not imply closedness. Consequently, Morse sets need not be closed and may be adjacent to one another.
By allowing adjacent Morse sets
we can detect finer Morse decompositions. We still can disconnect them
by modifying slightly the topology of the space without changing the topology of the Morse sets.

The organization of the paper is as follows.
In Section~\ref{sec:prelim} we recall preliminary material and
notation needed in the paper.
In Section~\ref{sec:comb-dyn} we introduce the concept of a  combinatorial dynamical system, define solutions and invariant sets of
a combinatorial dynamical system and present two
methods for constructing combinatorial dynamical systems
from data.
In Section~\ref{sec:IIS-MorseDecomp}
we define the concepts of isolating neighborhood, isolated invariant set and Morse decomposition of a combinatorial dynamical system.
In Section~\ref{sec:persistence} we define homological persistence of Morse
decompositions in the setting of combinatorial dynamical systems.
In Section~\ref{sec:geom-int} we discuss computational aspects of the theory
and provide a geometric interpretation of the Alexandrov topology
of subsets of a simplicial complexes.
In Section~\ref{sec:examples} we present two numerical examples.

\section{Preliminaries}
\label{sec:prelim}

In this section we recall some definitions and results needed in the paper
and establish some notations.

\subsection{Finite topological spaces}
\label{ssec:fts}
We recall that a {\em topology} on a set $X$ is a family $\cT$ of subsets
of $X$ which is closed under finite intersection and
arbitrary union and satisfies $\emptyset,X\in\cT$.
The sets in $\cT$ are called {\em open}.
The {\em interior} of $A$, denoted $\inte A$, is the union of open subsets of $A$.
A set $A\subseteq X$ is {\em closed} if $X\setminus A$ is open.
The {\em closure} of $A$, denoted $\cl A$, is the intersection of all closed supersets of $A$.
Topology $\cT$ is $T_2$
or Hausdorff if for any $x,y\in X$ where $x\neq y$, there exist disjoint sets
$U,V\in\cT$ such that $x\in U$ and $y\in V$.
It is $T_0$ or Kolmogorov if for any $x,y\in X$ such that $x\neq y$ there exists a $U\in\cT$ containing precisely one of $x$ and $y$.

A {\em topological space} is a pair $(X,\cT)$ where $\cT$ is a topology on $X$. It is a {\em finite topological space} if $X$ is finite.
Finite topological spaces differ from general topological spaces because the only Hausdorff topology
on a finite topological space $X$
is the discrete topology consisting of all subsets of $X$.

Given two topological spaces $(X,\cT)$ and $(X',\cT')$ we say that a map $f:(X,\cT)\to (X',\cT')$ is {\em continuous}
if  $U\in\cT'$ implies $f^{-1}(U)\in\cT$.

A remarkable feature of finite topological spaces is the following
theorem.
\begin{thm}
  \label{thm:alexandroff}
  (P. Alexandrov, \cite{Al1937})
  For a partial order $\leq$ on a finite set $X$, there is a $T_0$ topology $\cT_\leq$ on $X$
  whose open sets are upper sets with respect to $\leq$
  that is sets $T\subset X$ such that $x\in T$, $x\leq y$ implies
  $y\in T$.
  For a $T_0$ topology $\cT$ on $X$, there is a partial order $\leq_\cT$ where $x\leq_\cT y$
  if and only if $x$ is in the closure of $y$ with respect to $\cT$. The correspondences $\cT\mapsto\leq_\cT$
  and $\leq\mapsto\cT_\leq$ are mutually inverse. They transform continuous maps into an order-preserving maps and vice versa.
\end{thm}

The space $X$ is $\cT$-{\em disconnected} if there exist disjoint, non-empty sets $U,V\in\cT$ such that $X = U\cup V$.
The space $X$ is $\cT$-{\em connected}
if it is not $\cT$-disconnected. A subset $A\subset X$ is $\cT$-{\em connected} if it is connected as a space with induced topology $\cT_A$.
The connected component of $x\in X$, denoted $[x]_{\cT}$, is the union of
all connected subsets of $X$ containing $x$. Note, that $[x]_{\cT}$ is a connected set and
$\setof{[x]_{\cT}\mid x\in X}$ is a partition of $X$.

A {\em fence} in a poset $X$ is a sequence $x_0,x_1,\ldots,x_n$
of points in $X$ such that any two consecutive points are comparable.
$X$ is {\em order-connected} if for any two points $x,y\in X$ there exists a fence starting in $x$ and ending
in $y$.

\begin{prop}
\label{prop:connected-finite}
(\cite[Proposition 1.2.4]{Ba2011})
Let $(X,\cT)$ be a finite topological space. Then, the following conditions are equivalent:
\begin{itemize}
   \item[(i)] $X$ is a connected topological space.
   \item[(ii)] $X$ is order-connected with the preorder $\leq_\cT$.
   \item[(iii)] $X$ is a path-connected topological space.
\end{itemize}
\qed
\end{prop}

\subsection{Simplicial complexes as finite topological spaces}
Let $K$ be a finite simplicial complex, either a geometric
simplicial complex in $\RR^d$
(see \cite[Section 1.2]{Munkres1984}) or an abstract
simplicial complex (see \cite[Section 1.3]{Munkres1984}).
We consider $K$ as a poset $( K,\preceq)$ with $\sigma\preceq\tau$ if and only if $\sigma$ is a face of $\tau$
(also phrased $\tau$ is a {\em coface} of $\sigma$).
We define an {\em interval} between $\sigma$ and $\sigma'$, denoted by $[\sigma, \sigma']$, as a set $\{ \tau \in K \mid \sigma \preceq \tau \preceq \sigma'\}$.
The poset structure of $K$ provides, via Theorem~\ref{thm:alexandroff} (Alexandrov Theorem), a $T_0$ topology on $K$. We denote it $\cT_ K$ and
we refer to $\cT_ K$ as the {\em Alexandrov topology} of $K$. We note that $\cT_K$ is non-Hausdorff unless $ K$ consists of vertices only.
However, $\cT_ K$ is always $T_0$.

It is easy to see that a set $ A\subset  K$ is closed in the Alexandrov topology
if and only if all faces of any element of $ A$ are also in $ A$. Hence, the closure of $ A$
is the collection of all faces of elements in $ A$.

The non-Hausdorff topology $\cT_K$ of a simplicial complex $ K$  should not be confused
with the topology of the {\em polytope} $|K|$ of $ K$.
In the case of a geometric simplicial complex, the
polytope $|K|$ is just the union of all  simplices in $K$.
In the case of an abstract simplicial complex,
the polytope $|K|$ is defined up to a homeomorphism
as the polytope of a geometric realization of $K$
(see \cite[Sec. 1.2,1.3]{Munkres1984}).
Polytope $|K|$ is a subset of the Euclidean space with metric topology,
therefore its topology is Hausdorff.

An {\em open cell} $\cell{\sigma}$ associated with a simplex
$\sigma\in  K$ is the set of points $x$ in the polytope $| K|$
whose  barycentric coordinates $t_v(x)$ are strictly positive for every vertex $v\in\sigma$.
The {\em solid} of a set of simplices $ A\subset  K$ is $|A|:=\bigcup\setof{\cell{\sigma}\mid\sigma\in  A}$.
Note that $ A$ is a subcomplex of $K$ if and only if $A$ is closed in the Alexandrov topology of $K$ and then the solid of $A$ coincides with  the polytope of $A$.
This is why we use $|\cdot|$ to denote both solids and polytopes.
It is not difficult to verify that $ A\subset  K$ is open (respectively closed) in the Alexandrov topology if and only if
its solid is open (respectively closed) in $| K|$.

In the case of a geometric simplicial complex
we say that a set of simplices $A\subset K$ is {\em convex} if its solid $|A|$ is a convex set in $\RR^d$.
If the geometric simplicial complex $K$ is convex
and  $A\subset K$
we define the convex hull of $A$ as the intersection of all convex supersets of $A$ in $K$. We denote the
convex hull of $A$ by $\co A$.

\subsection{Multivalued maps and multivalued dynamics}
\label{ssec:mvd}
Recall that a {\em multivalued map} $F:X\mto Y$ is a map which assigns to every
point $x\in X$ a non-empty set $F(x)\subset Y$.
Given $A\subset X$, the {\em image} of $A$ under $F$ is
\[
      F(A):=\bigcup_{x\in A} F(x).
\]

For the sake of this paper we define a {\em multivalued dynamical system} in
a topological space $X$ as a multivalued map $F:X\times \NN \mto X$
such that
\begin{equation}
\label{eq:mds}
    F(F(x,m),n)=F(x,n+m).
\end{equation}
Typically, one also assumes that $F$ is continuous in some sense but we do not need such an assumption in this paper.

Let $F$ be a multivalued dynamical system. Consider the multivalued map $F^n: X\mto X$ given by
$F^n(x):=F(x,n)$. We call $F^1$ the {\em  generator} of
the dynamical system $F$.  It follows from \eqref{eq:mds} that the multivalued dynamical system $F$ is uniquely determined by the generator. Thus, it is natural to
identify a multivalued dynamical system with its generator.
In particular, we consider any multivalued map $F:X\mto X$
as a multivalued dynamical system $F:X\times \NN\mto X$
defined recursively by
\begin{eqnarray*}
   F(x,1)&:=&F(x)\\
   F(x,n+1)&:=&F(F(x,n)).
\end{eqnarray*}

\section{Combinatorial dynamics}
\label{sec:comb-dyn}
In this section we introduce the concept of a  combinatorial dynamical system and define solutions and invariant sets of
a combinatorial dynamical system. We also present two
cases for constructing combinatorial dynamical systems
from data.

\subsection{Combinatorial dynamical systems}
The central object of interest of this paper is given by the following definition.
\begin{defn}
\label{defn:cds}
{\em
   By a {\em combinatorial dynamical system} we mean
    a multivalued dynamical system generated by a multivalued map
   $F:X\mto X$ from a finite topological space $X$ to itself.
}
\end{defn}
In the sequel we identify the combinatorial dynamical system with its generator.
Although in this paper we restrict the considered examples to the case
of combinatorial dynamical systems generated by multivalued maps on
the collection of simplexes of a simplicial complex with its Alexandrov topology, the theoretical results apply to the general setting of finite topological spaces.
The general setting of finite topological space is useful, because there are methods to represent combinatorially subsets of $\RR^d$ other than the polytope of a simplicial complex, for instance
a cubical complex or a more general cellular complex.
All these cases lead to a finite topological space.
As we already mentioned in Section~\ref{ssec:mvd},
we do not require any continuity conditions on $F$.
Surprisingly, although such conditions are needed to define the Conley index (see \cite{MrWa2019}),
they are not needed to define  the isolating
neighborhood and Morse decomposition.

\subsection{Solutions and invariant sets}
A {\em solution} of $F$ in $A\subset K$ is a partial map $\rho:\ZZ\pto A$ whose {\em domain},
denoted $\dom \rho$,  is either the set of all integers
or a finite interval of integers and for any $i,i+1\in\dom\rho$ the inclusion $\rho(i+1)\in F(\rho(i))$ holds.
The solution $\rho$ is {\em full} if $\dom \rho=\ZZ$, otherwise it is {\em partial}.
In the latter case, if $\dom\rho=\ZZ\cap [m,n]$
for some $m,n \in \ZZ$, then $\rho(m)$ and $\rho(n)$ are called respectively
the left and right {\em endpoint} of $\rho$.
The solution {\em passes} through $\sigma\in K$ if $\sigma=\rho(i)$ for some $i\in\dom\rho$.
The set $A$ is {\em invariant} if for every $\sigma\in A$
there exists a full solution in $A$ passing through $\sigma$.

\subsection{A combinatorial dynamical system from a sampled map}
Assume $K$ is a convex simplicial complex in $\RR^d$ and consider a map $f:|K|\to |K|$ on the polytope of $ K$.
Moreover, assume we know only a noisy sample of $f$
that is a non-empty collection of pairs $\setof{(x_i,y_i)}_{i=1,n}$
satisfying $x_i,y_i\in | K|$ and $y_i$ equals $f(x_i)$ perturbed by some noise. Our goal is to investigate
the dynamical system generated by $f$ on $| K|$ by studying
a multivalued dynamical system induced on the finite collection
of simplexes of $ K$ by a multivalued map $ F: K\mto K$  constructed from the sample. In order to construct $ F$ recall that a {\em maximal simplex} or {\em toplex} in $K$ is a simplex
which is not a proper face of another simplex in $K$.
Denote by $ K_{top}$ the family of all toplexes in $ K$
and assume that each toplex is $d$-dimensional.
For toplexes $\tau,\tau'$ let $n_{\tau,\tau'}$ denote the number of pairs $(x_i,y_i)$ such that $x_i\in \cl |\tau|$ and $y_i\in \cl |\tau'|$.
Set $n_{max}:=\max\setof{n_{\tau,\tau'}\mid \tau,\tau'\in K_{top}}$
and let
\[
   \bar{n}_{\tau,\tau'}:=\frac{n_{\tau,\tau'}}{n_{max}}
\]
denote the relative frequency.
We first assign to each toplex $\tau$ the family
\[
A_{\mu,\tau}:=\setof{\tau'\in  K_{top}\mid
    \bar{n}_{\tau,\tau'} \geq\mu}
\]
that is the collection of toplexes $\tau'$ for which the relative frequency $\bar{n}_{\tau,\tau'}$ exceeds the threshold $\mu$.
Note that when the map $f$ is strongly expanding and the number of sample
points in $\tau$ is small, it may happen that
some or even all toplexes in $A_{\mu,\tau}$ are disjoint.
This is in contrast to the fact that
a continuous map sends a connected set to a connected set.
We remedy the problem by replacing $A_{\mu,\tau}$ with $\co A_{\mu,\tau}$, the convex hull of $A_{\mu,\tau}$. Next, we extend the definition to
a multivalued map $F_\mu: K\mto  K$ for an arbitrary
simplex $\sigma\in K$ (not necessarily a toplex) by setting
\[
    F_\mu(\sigma):=
    \co\bigcup\setof{A_{\mu,\tau}\mid \text{ $\sigma$ is a face of a toplex $\tau$} }.
\]

The multivalued map $F:=F_\mu$ is an example of the generator of a combinatorial dynamical system on the set of simplexes of the simplicial complex
$ K$.

\begin{figure*}[ht!]
  \begin{center}
    \includegraphics[width=0.3\textwidth]{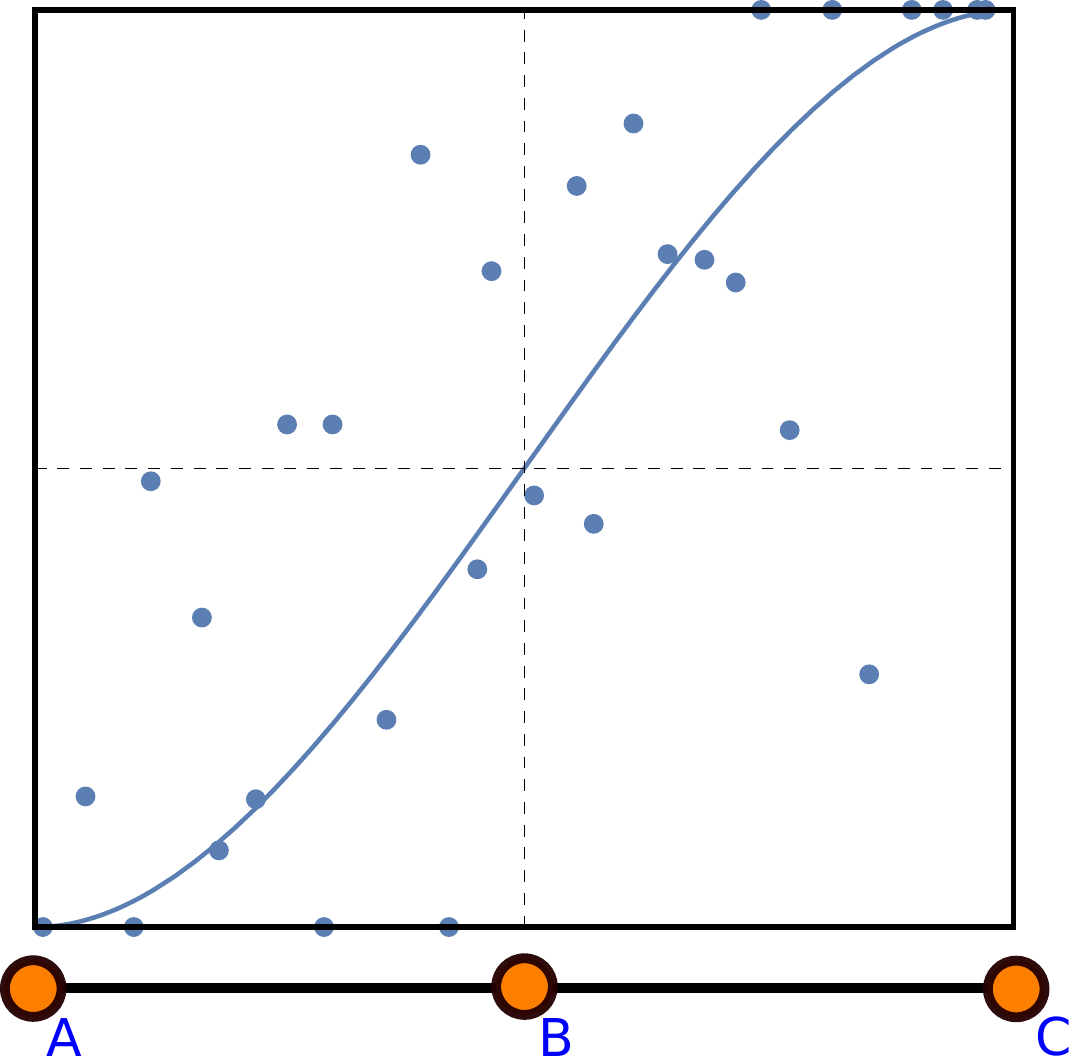}\hspace*{10pt}
    \raisebox{0.6\height}{\includegraphics[width=0.3\textwidth]{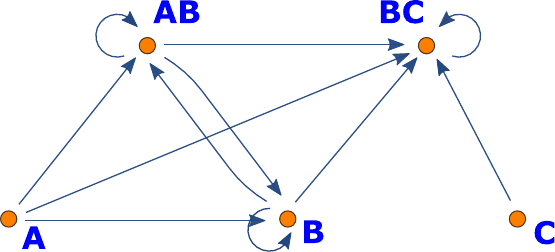}}\hspace*{10pt}
    \raisebox{0.6\height}{\includegraphics[width=0.3\textwidth]{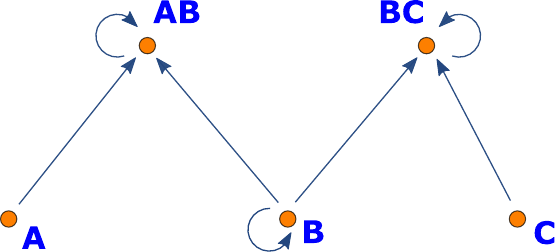}}
  \end{center}
  \caption{
  Bottom left: A simplicial complex in $\RR$ whose polytope
  is the interval $[0,1]$.
  Left: A map $f:[0,1]\ni x \mapsto 3x^2-2x^3\in [0,1]$ and a sample of $f$ with large Gaussian noise.
  Middle and right: The constructed combinatorial dynamical system $F_\mu$ with threshold $\mu=0.3$ (middle) and $\mu=0.4$ (right).
}
  \label{fig:s-map}
\end{figure*}

\subsection{A digraph interpretation of a combinatorial dynamical system}
A combinatorial dynamical system $F$ may be viewed as a digraph $G_F$
whose vertices are simplices in $ K$ with a directed edge from $\sigma$ to $\tau$ if and only if
$\tau\in F(\sigma)$.
An example is presented in Figure~\ref{fig:s-map}.
The polytope is the interval $[0,1]\subset\RR$.
The simplicial complex $K$
(see Figure~\ref{fig:s-map}(bottom left)) consists of two toplexes
$\bA\bB$ and $\bB\bC$ and three vertices $\bA$, $\bB$, $\bC$ where $\bA$, $\bB$,
$\bC$ are points in $\RR$ with coordinates $0$, $\frac{1}{2}$ and $1$ respectively.
A map $f:[0,1]\to [0,1]$ and a noisy sample
of this map are presented in Figure~\ref{fig:s-map}(left).
The relative freqencies are $\bar{n}_{AB,AB}=\frac{11}{12}$, $\bar{n}_{AB,BC}=\frac{4}{12}$, $\bar{n}_{BC,AB}=\frac{3}{12}$ and 
$\bar{n}_{BC,BC}=1$.
Figure~\ref{fig:s-map}(middle) and Figure~\ref{fig:s-map}(right) show digraph presentations
of two combinatorial dynamical systems on $ K$ respectively for thresholds $\mu=0.3$ and $\mu=0.4$.
In order to explain the presence of
the loop at vertex $\bB$ notice that $\bB$ is a face of two toplexes: $\bA\bB$
and $\bB\bC$.  For thresholds $\mu\in\{0.3,0.4\}$ we have
$\bA\bB\in A_{\mu,\bA\bB}$ and $\bB\bC\in A_{\mu,\bB\bC}$.
Therefore, $F_\mu(\bB)=\co\{\bA\bB,\bB\bC\}=\{\bA\bB,\bB,\bB\bC\}$.
But, an analogous computation for vertex $\bA$ gives
$F_\mu(\bA)=\co\{\bA\bB\}=\{\bA\bB\}$, which means that there is no loop at vertex $\bA$.
Similarly, we see that
there is no loop at vertex $\bC$.

The digraph interpretation of a combinatorial dynamical system means that
some concepts in dynamics may be translated into concepts in digraphs and vice versa.
In this translation a solution to $F$ in $A\subset K$ corresponds to a walk in $G_F$ through vertices in $A$
and the set $A$ is invariant if every vertex in $A$ is incident to a bi-infinite walk in $G_F$ through
vertices in $A$.
For instance, in Figure~\ref{fig:s-map}(middle),
the set $\{\bA\bB,\bB,\bB\bC\}$ is invariant.
Actually, all its subsets are also invariant
because of the presence of loops at $\bA\bB$, $\bB$ and
$\bB\bC$. The same comment applies to Figure~\ref{fig:s-map}(right).

We emphasize that, despite the convenience of the language of digraphs,  the combinatorial dynamical system $F$ is more than just the digraph $G_F$, because the collection of simplices $K$, that is the set of vertices of $G_F$,
is a $T_0$ topological space. In particular, the concept of isolating neighborhood which we define in Section~\ref{ssec:iso-inv-set},
cannot be formulated in the language of digraphs only.

%
\subsection{A combinatorial dynamical system from a sampled vector field}
When the dynamics which is sampled constitutes a flow, that is when time is continuous as in the case of a differential equation, the sampled data
often consists of a cloud of points with a vector attached to every point.
In this case the construction of combinatorial dynamical system is done in two
steps. In the first step the cloud of vectors is transformed into a combinatorial vector field
in the sense of Forman \cite{Fo98a,Fo98} or its generalized version of combinatorial multivector field
\cite{Mr2016}.
We discuss one of the possible algorithms for the first step in Section~\ref{sec:exlv}.
In the second step, the combinatorial multivector field is transformed into a combinatorial dynamical system.
In order to explain the second step, we introduce some definitions.
Let $K$ be a simplicial complex.
We say that $A\subset  K$ is {\em orderly convex} if for any $\sigma_1,\sigma_2 \in A$ and $\tau\in  K$ the relations $\sigma_1\preceq \tau$ and $\tau\preceq \sigma_2$ imply $\tau\in A$.
We define a {\em multivector} as an orderly convex subset of $ K$ and
a {\em combinatorial multivector field} on $ K$ ({\em combinatorial multivector field} in short) as a partition $\cV$ of $ K$ into multivectors.
Note that this definition of a combinatorial multivector field is less restrictive than the one in \cite{Mr2016}.
Both definitions encompass the combinatorial vector field of Forman as a special case.
The definition of combinatorial multivector field in \cite{Mr2016} additionally requires that multivectors have a unique maximal element. This is not needed here.

Given a combinatorial multivector field $\cV$, we denote by $[\sigma]_\cV$ the unique $V$ in $\cV$ such that $\sigma\in V$.
We associate with $\cV$ a combinatorial dynamical system $F_\cV: K\mto  K$ given by
$
   F_\cV(\sigma):=\cl \sigma \cup [\sigma]_\cV.
$
Note that $F_\cV$ in general admits more solutions than $\Pi_\cV$ defined in \cite[Section 5.4]{Mr2016}.
In particular, each $\sigma\in  K$ is a {\em fixed point} of $F_\cV$, that is, $\sigma\in F_\cV(\sigma)$.
This may look like a drawback but actually it simplifies the theory and allows detecting and eliminating spurious fixed points
by the triviality of their Conley index \cite{KuLiMr2018}.

\begin{figure*}[ht!]
  \begin{center}
    \raisebox{0.02\height}{\includegraphics[width=0.23\textwidth]{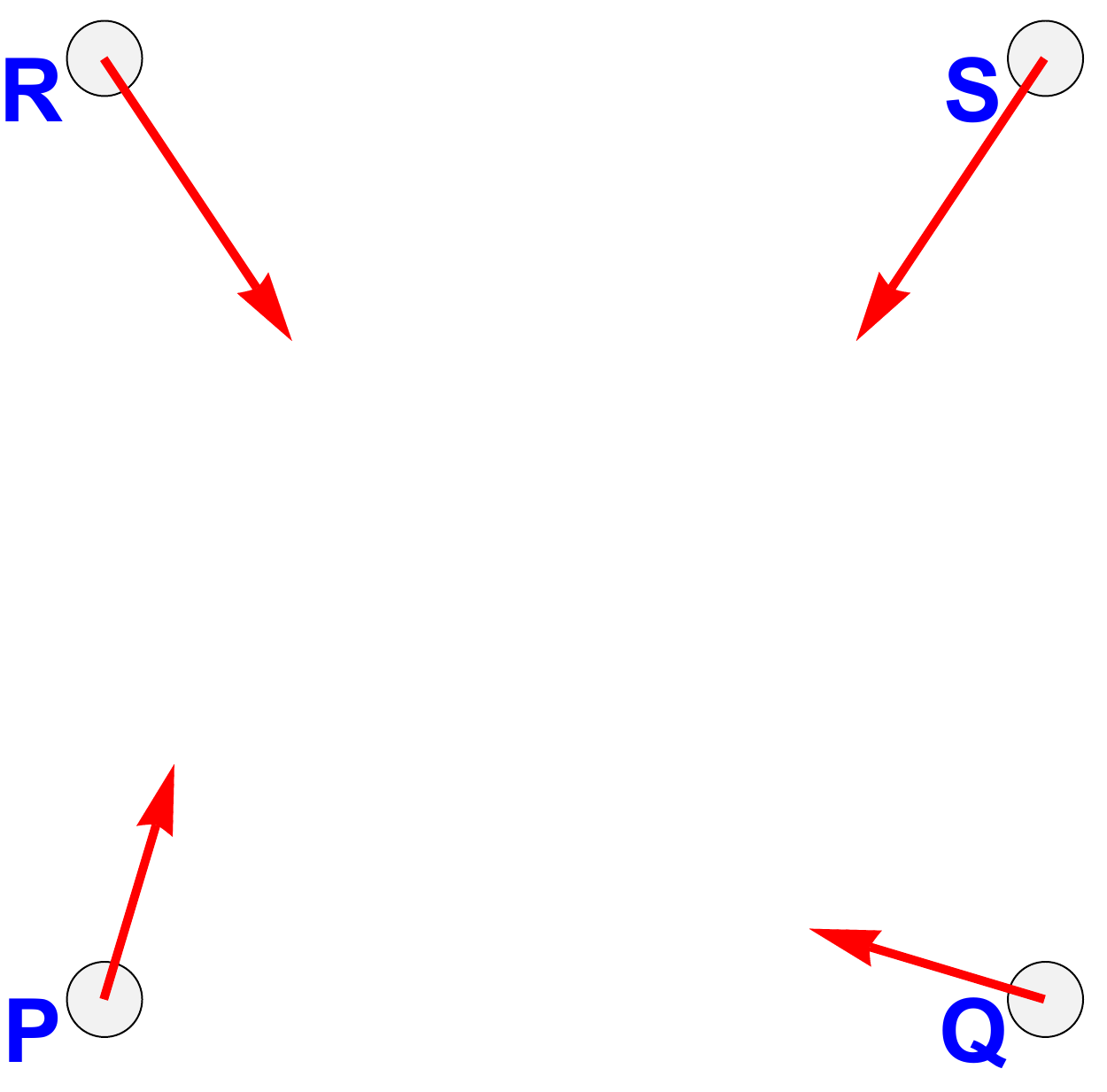}}\hspace*{15pt}
    \raisebox{0.02\height}{\includegraphics[width=0.23\textwidth]{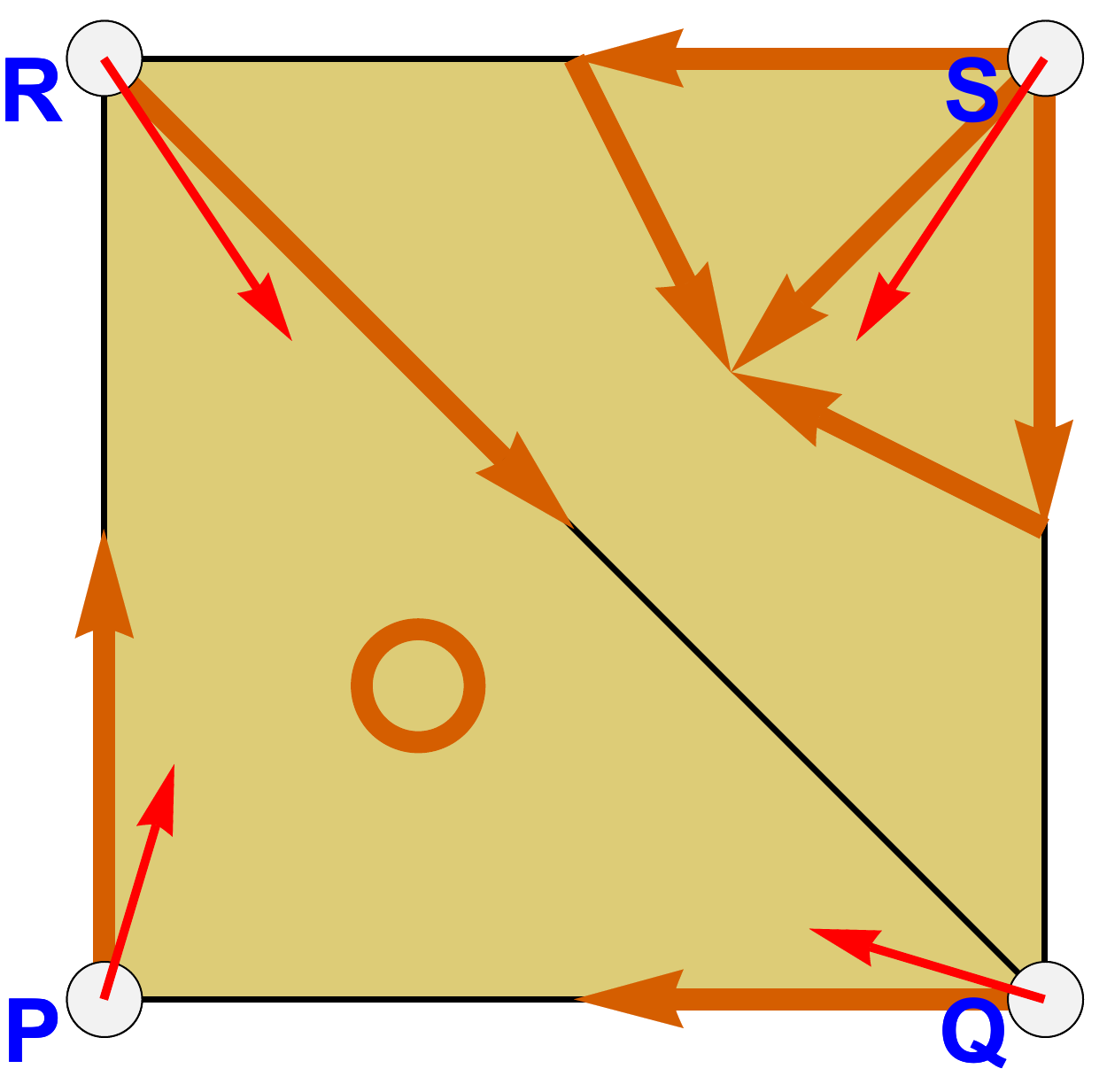}}\hspace*{10pt}
    \includegraphics[width=0.4\textwidth]{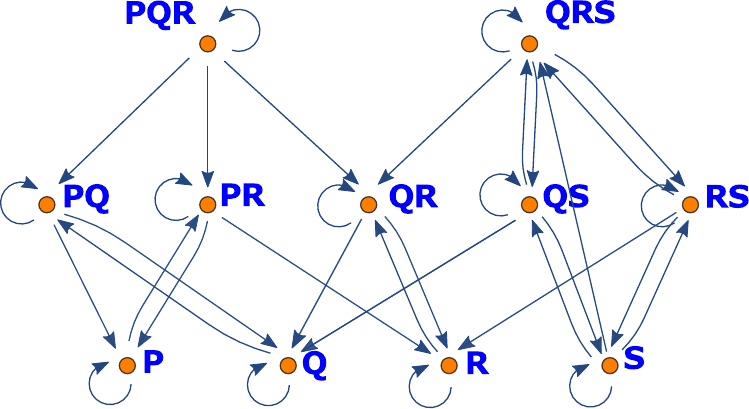}
  \end{center}
  \caption{Left: A cloud of vectors.
  Middle: A possible combinatorial multivector field representation of the cloud of vectors.
    Right: The associated combinatorial dynamical system represented as a digraph.}
  \label{fig:fv-map}
\end{figure*}

Figure~\ref{fig:fv-map}(left) presents a toy example of a cloud of vectors.
It consists of four vectors marked red at four points $\bP$, $\bQ$, $\bR$, $\bS$.
One of possible geometric simplicial complexes with vertices
at points $\bP$, $\bQ$, $\bR$, $\bS$ is the simplicial complex $K$ consisting of
triangles $\bP\bQ\bR$, $\bQ\bR\bS$ and its faces.
A possible multivector field $\cV$ on $ K$ constructed from the cloud of vectors consists of multivectors $\{\bP,\bP\bR\}$, $\{\bR,\bQ\bR\}$, $\{\bQ,\bP\bQ\}$, $\{\bP\bQ\bR\}$,$\{\bS,\bR\bS,\bQ\bS,\bQ\bR\bS\}$.
It is indicated in Figure~\ref{fig:fv-map}(middle) by orange
arrows between centers of mass of simplices.
Note that in order to keep the figure legible, only arrows in the direction
increasing the dimension are marked.
The singleton  $\{\bP\bQ\bR\}$ is marked with an orange circle.
The associated combinatorial dynamical system $F_\cV$ presented as a digraph is in Figure~\ref{fig:fv-map}(right).
Note that in general  $K$ and $\cV$ are  not uniquely determined by the cloud of vectors.
One possible method for constructing combinatorial multivector fields from a cloud of vectors is discussed
in Section~\ref{sec:exlv}.

\section{Isolated invariant sets and Morse decompositions}
\label{sec:IIS-MorseDecomp}

In this section we consider a combinatorial dynamical system $F:X\mto X$ and define for $F$ the concepts of isolating neighborhood, isolated invariant set and Morse decomposition.

\subsection{Isolated invariant sets}
\label{ssec:iso-inv-set}
The closed set $ N\subset K$ is an {\em isolating neighborhood} for an invariant set $S\subset K$
if $S$ is contained in  $N$ and any partial solution in $ N$ with endpoints in $S$
has all values in $S$. If such an isolating neighborhood for $S$ exists, we say that $S$
is an {\em isolated invariant set}.
We emphasize that, unlike the classical theory,
the same set $N$ may be an isolating neighborhood
for more than one isolated invariant set.

The invariant set $\{\bA\bB\}$ in Figure~\ref{fig:s-map}(middle)
is not an isolated invariant set, because for any closed set $N$ containing $\{\bA\bB\}$  the partial solution $(\bA\bB,\bB,\bA\bB)$ is contained in $N$ and has endpoints in $\{\bA\bB\}$.
The invariant sets $\{\bB\bC\}$ and $\{\bA\bB,\bB\}$ are both isolated invariant sets and
$\{\bA,\bA\bB,\bB,\bB\bC,\bC\}$ is an isolating neighborhood for both.

Since we have a loop at every vertex of the digraph of
the combinatorial dynamical system in Figure~\ref{fig:fv-map}(right), every set is invariant. In particular, every singleton
is invariant. However, the only singleton which is an isolated
invariant set is $\{\bP\bQ\bR\}$. Its isolating neighborhood
is $\cl\{\bP\bQ\bR\}=\{\bP,\bQ,\bR,\bP\bQ,\bP\bR,\bQ\bR,\bP\bQ\bR\}$.
Another isolated invariant set with the same isolating neighborhood is $\{\bP\bQ,\bP\bR,\bQ\bR\}$.

The {\em maximal invariant set} of $F$, denoted $S(F)$, is the set of all
simplices $\sigma\in K$ such that there exists
a full solution of $F$ in $K$ passing through $\sigma$.
 It is straightforward to observe that $S(F)$ is invariant
and $K$ is an isolating neighborhood for $S(F)$. Therefore, $S(F)$ is an isolated invariant set.
Note that the maximal invariant set $S(F_\cV)$ for a combinatorial multivector field $\cV$ is always the whole $K$,
because for each $\sigma\in K$ we have $\sigma\in\cl\sigma\subset F_\cV(\sigma)$.
This is visible in Figure~\ref{fig:fv-map}(right) as a loop at every vertex.
In contrast, $\bA$ does not belong to the maximal invariant set in Figure~\ref{fig:s-map}(right).

\subsection{Morse decompositions}
A {\em connection} from an isolated invariant set $S_1$ to an isolated invariant set $S_2$
is a partial solution with left endpoint in $S_1$ and right endpoint in $S_2$.
A family $\cM$ consisting of mutually disjoint, non-empty isolated invariant subsets of an isolated invariant set~$S$
is a {\em Morse decomposition} of~$S$ if $\cM$ admits a partial order $\leq$ such that
any connection between elements in~$\cM$ either has all values
in a single element of $\cM$ or it originates in $M\in\cM$ and terminates in $M'$ such that $M>M'$.
If $S$ is not mentioned explicitly, we mean a Morse decomposition of the maximal invariant set $S(F)$.
The elements of $\cM$ are called {\em Morse sets}.
Although the definitions of isolated invariant set and Morse decomposition require topology, there is an important
case when they correspond to purely graph-theoretic concepts.
An isolated invariant set is {\em minimal} if it admits no non-trivial Morse decomposition that is no Morse decomposition
consisting of more than one Morse set. A Morse decomposition is {\em minimal} if each of its Morse sets is minimal.
The following theorem shows that the minimal Morse decomposition of $F$, denoted as $\cM(F)$, is unique and
consists of the strongly connected components of $G_F$.
\begin{thm}
\label{thm:ssc-as-morse}
The family of all strongly connected components of $G_F$ is the unique minimal Morse decomposition of $S(F)$.
\end{thm}
\begin{proof}
   Let $\cS$ be the family of all strongly connected components of $G_F$.
   We will show that $K$ is an isolating neighborhood for any $S\in\cS$.
   Obviously, $K$, as the whole space  is  closed. Therefore,  $S\subset K$.
   Moreover, any partial solution with endpoints in $S$ must have all values in $S$, because $S$ is a strongly connected
   component of $G_F$.
   Hence, each $S\in\cS$ is an isolated invariant set. Clearly, it is a minimal isolated invariant set.
   For $S_1,S_2\in\cS$ we write $S_1\geq S_2$ if the there exists a connection from $S_1$ to $S_2$.
   Since $\cS$ consists of strongly connected components, this defines a partial order on $\cS$.
   Let $\rho$ be a connection from $S_1$ to $S_2$ whose
values are not contained in a single element of $\cS$.
Then, $S_1>S_2$. This proves that $\cS$ is a Morse decomposition.
Obviously, a strongly connected component cannot have
a non-trivial Morse decomposition. Thus, $\cS$
is a minimal Morse decomposition.
   Assume that $\cS'$ is another minimal Morse decomposition and $S'\in\cS'$. We claim that $S'$ is strongly connected
as a subgraph of $G_F$.
Indeed, if not, then, according to what we already proved,
the strongly connected components of $S'$ would constitute
a non-trivial Morse decomposition of $S'$,
contradicting the assumption that $S'$ is a minimal invariant set.
Hence, $S'$ is contained in a Morse set $S\in\cS$.
By a symmetric argument every $S\in \cS$ is contained
in a Morse set  $S'\in\cS'$.
   This shows that $\cS=\cS'$ and proves the uniqueness.
   \qed
\end{proof}

\begin{figure*}[ht!]
  \begin{center}
    \includegraphics[width=0.25\textwidth]{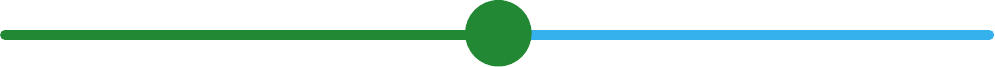} \hspace*{10pt}
    \includegraphics[width=0.25\textwidth]{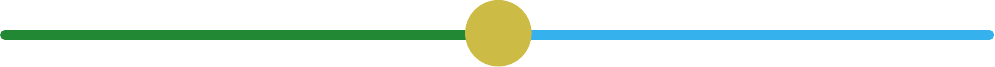} \hspace*{10pt}
    \includegraphics[width=0.25\textwidth]{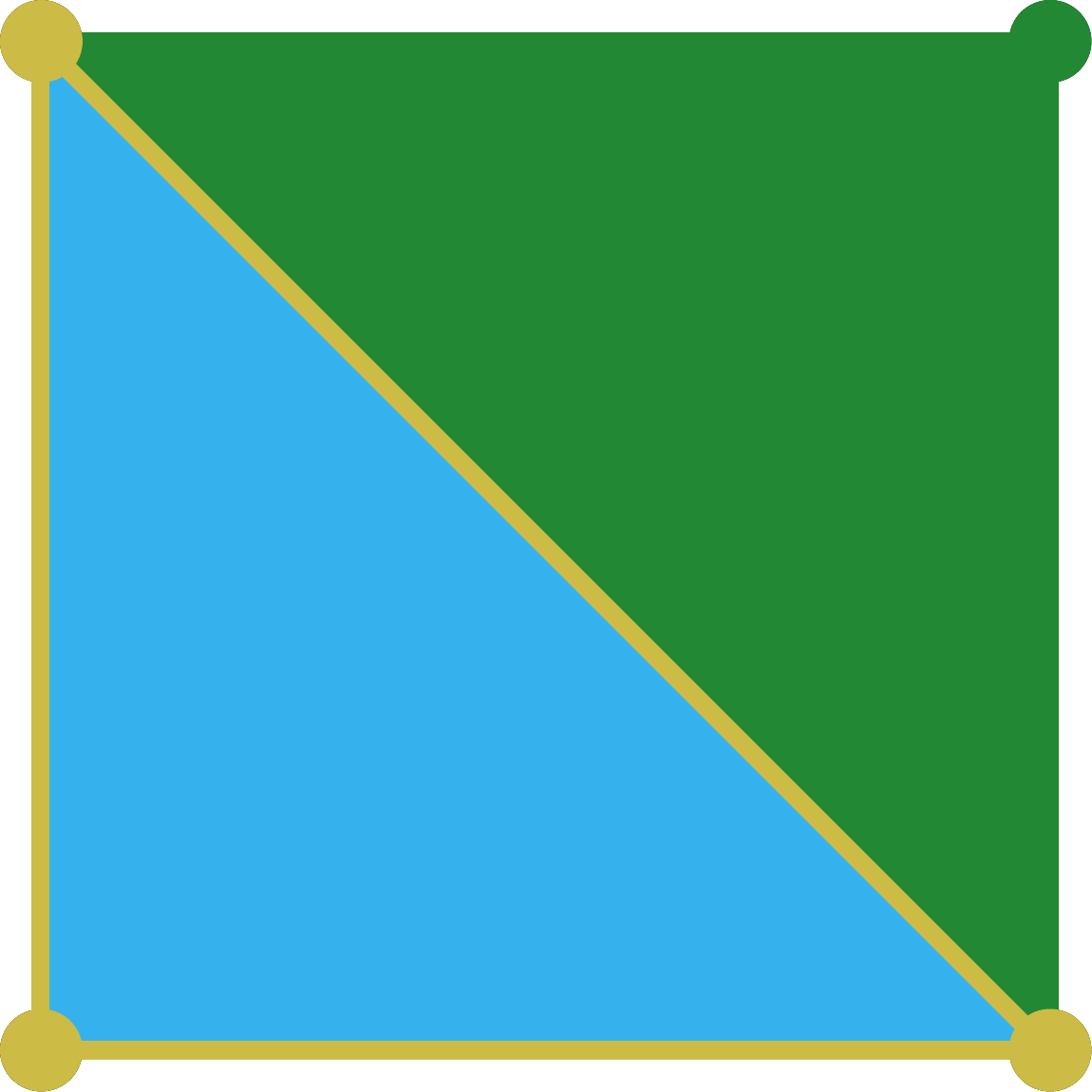}       \end{center}
  \caption{Left: The minimal Morse decomposition of the combinatorial
  dynamical system in Figure~\ref{fig:s-map}(middle) consisting of
  $\{\bA\bB,\bB\}$ (green) and $\{\bB\bC\}$ (blue). Middle: The minimal Morse decomposition of the combinatorial
  dynamical system in Figure~\ref{fig:s-map}(right) consisting of
  $\{\bA\bB\}$ (green), $\{\bB\}$ (yellow) and $\{\bB\bC\}$ (blue).
    Right: The minimal Morse decomposition of the combinatorial
  dynamical system in Figure~\ref{fig:fv-map}(right)
  consisting of $\{\bP,\bQ,\bR,\bP\bQ,\bP\bR,\bQ\bR\}$ (yellow),
$\{\bS,\bR\bS,\bQ\bS,\bQ\bR\bS\}$ (green) and $\{\bP\bQ\bR\}$ (blue).
  }
\label{fig:fv-map-md}
\end{figure*}

Consider the combinatorial dynamical system in Figure~\ref{fig:s-map}(middle).
Its minimal Morse decomposition
consists of two Morse sets: $\{\bA\bB,\bB\}$ and $\{\bB\bC\}$ with $\{\bA\bB,\bB\}>\{\bB\bC\}$.
The minimal Morse decomposition of the combinatorial dynamical system in Figure~\ref{fig:s-map}(right)
consists of three Morse sets: $\{\bA\bB\}$, $\{\bB\}$ and $\{\bB\bC\}$ with $\{\bB\}>\{\bB\bC\}$ and $\{\bB\}>\{\bA\bB\}$.
The minimal Morse decomposition of the combinatorial dynamical system in Figure~\ref{fig:fv-map}(right) consists of three isolated invariant sets:  $M_1:=\{\bP,\bQ,\bR,\bP\bQ,\bP\bR,\bQ\bR\}$,
$M_2:=\{\bS,\bR\bS,\bQ\bS,\bQ\bR\bS\}$ and $M_3:=\{\bP\bQ\bR\}$ with $M_3>M_1$ and $M_2>M_1$.
These minimal Morse decompositions are illustrated in Figure~\ref{fig:fv-map-md}.

\section{Persistence of Morse decompositions.}
\label{sec:persistence}
In this section we define homological persistence of Morse
decompositions in the setting of combinatorial dynamical systems.
\subsection{Disconnecting topology.}
In the case of a classical Morse decomposition consisting of more than one Morse set, the union of all Morse sets is always disconnected in the topology induced from the space.
This is because Morse sets are always disjoint and in this case also compact. In particular, the space between the Morse sets is filled with solutions connecting the Morse sets.
But, in finite topological spaces the Morse sets need not be closed
and solutions may jump directly from one Morse set to another Morse set. Consequently, the union of Morse sets
generally is not disconnected. Thus, we need a method to disconnect Morse sets. Fortunately, in this case we do not need space between
the Morse sets.
We achieve the separation  by purely topological methods. To explain this, we need the following terminology, notation and theorem.

Assume $\cA$ is a finite family of mutually disjoint non-empty sets and $\cT$ is a topology on $\bigcup \cA$. We say that $\cA$ is {\em disconnected} in $\cT$
if each set $A\in\cA$ is open in the topology $\cT$.

Given a family $\cA$ of subsets of a set $X$,
we use the notation $\cA^* := \{\bigcup \cA' : \cA' \subset \cA\}$ for the smallest family of sets in $X$,
containing $\cA$ and closed under summation.
If $\cB$ is another such family, we write $\cA\bcap\cB$ for the family of intersections
of every set in $\cA$ with every set in $\cB$.
We say that $\cA$ is {\em inscribed} in $\cB$ and write $\cA\sqsubset\cB$ if for every $A\in\cA$, there exists
a $B\in\cB$ such that $A\subset B$.

In order to shorten the notation
we will also write $\union{\cA}$ for the union $\bigcup\cA$
of all the sets in $\cA$.
Note that if $A\subset X$ and $\cT$ is a topology on $X$, then the topology induced by $\cT$ on $A$ is $A\bcap\cT:=\{A\}\bcap\cT$.

\begin{thm}
\label{thm:conn-comp-intersec}
Assume  $(X,\cT)$ is an arbitrary topological space and $\cA$ is a finite
family of mutually disjoint, non-empty subsets of $X$.
Then $\cT_\cA:=\left(\cA\bcap\cT\right)^*$ is a topology on $\union{\cA}$. Moreover,
\begin{itemize}
   \item[(i)] if $\cT$ is a $T_0$ topology, then so is $\cT_\cA$,
   \item[(ii)] for every $A\in\cA$, the topology induced on $A$ by $\cT$ coincides  with the topology induced on $A$ by $\cT_\cA$,
   \item[(iii)] the family $\cA$ is  $\cT_\cA$-disconnected,
   \item[(iv)] if additionally $\union{\cA}=X$ and each set in $\cA$ is $\cT$-connected, then the connected components
   with respect to  $\cT_\cA$ coincide with the sets in $\cA$.
\end{itemize}
\end{thm}
\begin{proof}
We will show that $\cA\bcap\cT$ is a basis (see \cite[Section 13]{Munkres1975})
for some topology  on $\union{\cA}$. Let $x \in \union{\cA}$.
There exists an $A \in \cA$ such that $x \in A$.
Hence, $x\in A = A \cap X \in \cA\bcap\cT$.
Assume that $x \in (A \cap U) \cap (B \cap V)$  for some $A, B \in \cA$ and $U, V \in \cT$.
Then $A = B$ and consequently $(A \cap U) \cap (B \cap V) = A \cap (U \cap V)  \in \cA\bcap\cT$.
This shows that $\cA\bcap\cT$ is indeed a basis. By \cite[13.1]{Munkres1975} it follows that $\cT_\cA$ is a topology.
Consider $x, y \in \union{\cA}$,  $x \neq y$.
Without loss of generality we can assume that there exists an open neighbourhood $U\in\cT$ of $y$ such that $x \not\in U$.
Let $A \in \cA$ be such that $y \in A$.
Then $y \in U \cap A$ and $x \not\in U \cap A$, hence (i) holds.
To prove (ii) we need to show that $A\bcap\cT=A\bcap\cT_\cA$. Obviously $A\bcap\cT\subset A\bcap\cT_\cA$.
To prove the opposite inclusion take a $V\in  A\bcap\cT_\cA$.
This means that there is a $U \in \cT_\cA$ such that $V=A\cap U$.
Then $U=\bigcup_{i\in I}(A_i \cap U_i)$ for some $U_i \in \cT$ and $A_i \in \cA$.
Hence $V=A\cap U=A\cap\bigcup_{i\in I}( A_i \cap U_i)  = A\cap \bigcup_{i\in I_A} U_i  \in A\bcap\cT$ where $I_A = \{i\in I : A_i = A\}$ and (ii) is proved.
Property (iii) is obvious, because  $A \in \cA$ implies $A = A \cap X \in \cT_\cA$.
To prove (iv) assume $A$ is $\cT$-connected.
It follows from (ii) that  $A$ is $\cT_\cA$-connected.
Let $x \in A$. Then $A$ is contained in $[x]_{\cT_\cA}$, the $\cT_\cA$-connected component of $x$.
This means that $[x]_{\cT_\cA}=\bigcup \cA'$ for some $\cA' \subset \cA$.
Since every set in $\cA$ is open in $\cT_\cA$, the family $\cA'$ must contain precisely one element.
Consequently $A=[x]_{\cT_\cA}$ and (iv) holds.
\qed
\end{proof}

We note that given a Morse decomposition $\cM$
of a combinatorial dynamical system $F$ on a finite
simplicial complex $K$,  in general  the union $\union{\cM}$ is not a subcomplex of the simplicial complex $K$.
Therefore, we cannot take simplicial homology of $\union{\cM}$. Moreover, we are interested
in the special topology $\cT_\cM$ on $\union{\cM}$ where $\cT$ is the Alexandrov topology of $K$. The topology $\cT_\cM$ separates the Morse sets due to Theorem~\ref{thm:conn-comp-intersec}(iii).
Fortunately, the singular homology makes
sense for any topological space, in particular we can consider $H(\union{\cM},\cT_\cM)$.
In Section~\ref{sec:geom-int} we use McCord's Theorem \cite{McCord1966} to show that
$H(\union{\cM},\cT_\cM)$ may be computed as simplicial homology of a subcomplex of the barycentric subdivision of $K$.

\subsection{Persistence and zig-zag persistence of Morse decompositions}
Consider  two simplicial complexes $K$ and $K'$ with combinatorial dynamical systems $F$ on $K$ and $F'$ on $K'$
and a  map $f:K\to K'$ continuous with respect to
Alexandrov topologies $\cT$ on $K$ and $\cT'$ on $K'$.
By Theorem~\ref{thm:alexandroff} (Alexandrov Theorem) the map $f$ is continuous if and only if it preserves the face relation in $K$ and $K'$.
In particular, every simplicial map is continuous.

The following theorem lets us define homomorphisms in homology needed to set up
persistence of Morse decompositions.

\begin{thm}
\label{thm:f-map-on-Morse}
Let $\cM$ and $\cM'$ be Morse decompositions respectively for $F$ and $F'$.
Assume that continuous map $f:K \to K'$ {\em respects} $\cM$ and $\cM'$ that is $f(\cM)\sqsubset\cM'$ where $f(\cM):=\setof{f(M)\mid M\in\cM}$.
Then, the map $f_{\cM,\cM'}:(\union{\cM},\cT_\cM)\ni\sigma\mapsto f(\sigma)\in (\union{\cM'},\cT'_{\cM'})$ is well defined and continuous.
\end{thm}
\begin{proof}
   Let $\sigma\in\union{\cM}$. Then $\sigma\in M$ for some $M\in\cM$.
   Since $f$ respects $\cM$ and $\cM'$, there is an $M'\in\cM'$ such that $f(M)\subset M'$.
   It follows that $f(\sigma)\in \union{\cM'}$. Hence, $f_{\cM,\cM'}$ is well defined.
   Since $\cM'\bcap\cT'$ is a basis of topology $\cT'_{\cM'}$, in order to prove continuity it suffices to show that for any $M'\in\cM'$ and
   $T'\in\cT'$ the set $f_{\cM,\cM'}^{-1}(M'\cap T')$ is open in $\cT_\cM$. Let $\cM_{M'}:=\setof{M\in\cM\mid f(M)\subset M'}$.
   Then $f^{-1}(M')\cap\union{\cM}=\union{\cM_{M'}}$. By continuity of $f$ we have $f^{-1}(T')\in\cT$. Therefore,
$
  f_{\cM,\cM'}^{-1}(M'\cap T')=f^{-1}(M')\cap f^{-1}(T')\cap\union{\cM}=
  f^{-1}(T')\cap\union{\cM_{M'}}=
  \bigcup\{M\cap f^{-1}(T')\mid M\in\cM_{M'}\}\in(\cM\bcap\cT)^*=\cT_\cM,
$
which completes the proof.
\qed
\end{proof}

For a minimal Morse decomposition, denoted by $\cM(F)$, we have the following corollary.

\begin{cor}
\label{cor:f-map-on-Morse}
  The map $f_{\cM(F),\cM(F')}:(\union{\cM(F)},\cT_{\cM(F)}) \to (\union{\cM(F')},\cT'_{\cM(F')}) $ is continuous
  under the assumption that $f\circ F\subset F'\circ f$ that is $f(F(\sigma))\subset F'(f(\sigma))$ for any $\sigma\in K$.
\end{cor}
\begin{proof}
  By Theorem~\ref{thm:f-map-on-Morse}, it suffices to show that $f$ respects $\cM(F)$ and $\cM(F')$.
  Let $M\in \cM(F)$. By Theorem~\ref{thm:ssc-as-morse}, the Morse set $M$ is a strongly connected component
  of $G_F$. Let $\sigma,\tau\in M$ and let $\rho$ be a partial solution in $M$ with endpoints $\sigma$ and $\tau$.
  It follows from the assumption that $f\circ\rho$ is a solution in $f(M)$ with endpoints $f(\sigma)$ and $f(\tau)$.
  Since $\sigma,\tau\in M$ are arbitrary, the set $f(M)$ must be contained in one strongly connected component
  of $G_{F'}$, that is $f(M)\subset M'$ for some $M'\in\cM'$.
\qed
\end{proof}

Assume now that for $i=1,2,\ldots n$, we have a simplicial complex $K_i$ with Alexandrov topology $\cT^i$, a combinatorial dynamical system $F_i$ on $K_i$,
and a Morse decomposition $\cM_i$ of $F_i$.
Let  $\{f_i: K_i\to K_{i+1}\}_{i=1,n-1}$ be a sequence of continuous maps such that $f_i\circ F_i\subset F_{i+1}\circ f_i$
and $f_i(\cM_i)\sqsubset\cM_{i+1}$.
Note that by Corollary~\ref{cor:f-map-on-Morse} the latter condition may be dropped if $\cM_i=\cM(F_i)$.
It follows from Theorem~\ref{thm:f-map-on-Morse}
that the maps $\bar{f}_i:=(f_i)_{\cM_i,\cM_{i+1}}: (\union{\cM_i},\cT^i_{\cM_i}) \to (\union{\cM_{i+1}},\cT^{i+1}_{\cM_{i+1}})$ are continuous.
Thus, we have homomorphisms induced in singular homology
$H(\bar{f}_i):H(\union{\cM_i},\cT^i_{\cM_i}) \to H(\union{\cM_{i+1}},\cT^{i+1}_{\cM_{i+1}})$.
This yields a persistence module
\begin{equation}
\label{eq:persistence-module-simp}
\xymatrix{
  H(\union{\cM_1},\cT^1_{\cM_1}) \ar[r]^-{H( \bar{f}_1)} & H(\union{\cM_2},\cT^2_{\cM_2})\ar[r]^-{H( \bar{f}_2)} &  \dots \ar[r]^-{H( \bar{f}_{n-1})} & H(\union{\cM_n},\cT^n_{\cM_n}).
}
\end{equation}
We refer to the persistence diagram of this module as the {\em persistence diagram of Morse decompositions}.
We note that zig-zag persistence diagram of Morse decompositions may be obtained analogously
by replacing, whenever appropriate, inclusions  $f_i\circ F_i\subset F_{i+1} \circ f_i$ by $f_i\circ F_i\supseteq F_{i+1} \circ f_i$
and respectively $\bar{f}_i(\cM_i)\sqsubset\cM_{i+1}$ by $\cM_i\sqsupset\bar{f}_i(\cM_{i+1})$.


\subsection{Persistence in combinatorial multivector fields}
Let $\cV$ be a combinatorial multivector field on a simplicial complex $K$.
We say that $\cM$ is a Morse decomposition of $\cV$ if it is a Morse decomposition of the associated combinatorial dynamical system $F_\cV$.
We extend this terminology to minimal Morse decompositions.
We denote the minimal Morse decomposition of $\cV$ by $\cM(\cV):=\cM(F_\cV)$
and the topology of this Morse decomposition by $\cT_\cV:=\cT_{\cM(\cV)}$.

\begin{thm}
Morse decompositions of combinatorial multivector fields have the following properties.
\begin{itemize}
   \item[(i)] The minimal Morse decomposition of a combinatorial multivector field $\cV$ on $K$ is a partition of $K$. In particular, $\union{\cM}=K$.
   \item[(ii)] Given $\cW$, another combinatorial multivector field on $K$, the family $\cV\bcap \cW$ is a combinatorial multivector field on $K$. It is inscribed both in $\cV$ and $\cW$.
              Moreover, If $\cV\sqsubset \cW$, then $F_\cV\subset F_\cW$.
   \item[(iii)] If $\cV'$ is a combinatorial multivector field on a simplicial complex $K'$ and $f:K\to K'$ is continuous, then
                $
                f^*(\cV'):=\setof{f^{-1}(V')\mid V'\in\cV'},
                $
                called the {\em pullback} of $\cV'$, is a combinatorial multivector field on $K$.
   \item[(iv)] The maps
                $
                \kappa:=\id_{\cV\bcap f^*(\cV'),\cV}:(K,\cT_{\cV\bcap f^*(\cV')})\to (K,\cT_{\cV})
                $
                induced by identity and
                $
                \lambda:=f_{\cV\bcap f^*(\cV'),\cV'}:(K,\cT_{\cV\bcap f^*(\cV')})\to (K',\cT'_{\cV'})
                $
                induced by $f$ are continuous.
\end{itemize}
\end{thm}
\begin{proof}
  Note that by Theorem~\ref{thm:ssc-as-morse}, the Morse sets in the minimal Morse decomposition
  are the strongly connected components of $G_{F_\cV}$. Hence, to prove (i) it suffices to observe
  that every $\sigma\in K$ belongs to a strongly connected component. This is obvious
  because $\sigma\in\cl\sigma\subset F_\cV(\sigma)$ for any $\sigma\in K$.
  Thus, (i) is proved.
  Since the intersection of two orderly convex sets is easily seen to be orderly convex, each element of $\cV\bcap \cW$ is orderly convex.
  Obviously, $\cV\bcap \cW$ is a partition of $K$ and is inscribed in $\cV$ and $\cW$.
  Take $\sigma\in K$. Assumption $\cV\sqsubset \cW$ implies that $[\sigma]_\cV\subset[\sigma]_\cW$.
  It follows that $F_\cV(\sigma)=\cl\sigma\cup[\sigma]_\cV\subset \cl\sigma\cup[\sigma]_\cW=F_\cW(\sigma)$.
  Thus, (ii) is also proved.
  Obviously, $f^*(\cV')$ is a partition of $ K$.
  To show that for every $V'\in\cV'$ the set $f^{-1}(V')$ is orderly convex,
  take $\sigma,\sigma'\in f^{-1}(V')$
  and $\tau\in K$ such that $\sigma\preceq\tau\preceq\sigma'$. Then $f(\sigma),f(\sigma')\in V'$,
  $f(\sigma)\preceq f(\tau)\preceq f(\sigma')$, and since $V'$ is orderly convex, we get $f(\tau)\in V'$.
  It follows that $\tau\in f^{-1}(V')$ and $f^{-1}(V')$ is orderly convex. This proves (iii).
  To prove (iv), we verify that the maps $\kappa$ and $\lambda$  satisfy the assumption of
  Corollary~\ref{cor:f-map-on-Morse}.
  It follows from (ii) that
$
   \id\circ F_{\cV\bcap f^*(\cV')} = F_{\cV\bcap f^*(\cV')}\subset F_{\cV} = F_{\cV}\circ \id
$
which proves that $\kappa$ is continuous. Similarly, we get
$
   f\circ F_{\cV\bcap f^*(\cV')}\subset f\circ F_{f^*(\cV')}.
$
Thus, it suffices to prove that $f\circ F_{f^*(\cV')}\subset F_{\cV'}\circ f$. Indeed, for $\sigma\in K$ we get from the continuity of $f$
and the definition of $f^*(\cV')$ that
$
(f\circ F_{f^*(\cV')})(\sigma)=f(F_{f^*(\cV')}(\sigma))=f(\cl\sigma\cup[\sigma]_{f^*(\cV')})=
f(\cl\sigma)\cup f([\sigma]_{f^*(\cV')})\subset
\cl f(\sigma)\cup [f(\sigma)]_{\cV'}=F_{\cV'}(f(\sigma))=(F_{\cV'}\circ f)(\sigma).
$
\qed
\end{proof}

We use the diagram of continuous maps
$
\xymatrix{
 (K,\cT_\cV)  &  (K,\cT_{\cV\bcap f^*(\cV')}) \ar[l]_-{\kappa} \ar[r]^-{\lambda} & (K',\cT_{\cV'}),
}
$
referred to as the {\em comparison diagram} of combinatorial multivector fields $\cV$ and $\cV'$, to define the persistence of Morse decompositions
for combinatorial multivector fields.
To this end, assume that, for $i=1,2,\ldots n$, we have a combinatorial multivector field\ $\cV_i$ on a simplicial complex $K_i$.
Moreover, assume that we have a sequence of continuous maps $f_i:K_i\to K_{i+1}$.
Putting together the comparison diagrams of $\cV_i$ and $\cV_{i+1}$ and applying the singular homology functor
we obtain the following zig-zag persistence module
\begin{multline}
\label{eq:fi-zigzag-diagram1}
\xymatrix{
  H(K_1,\cT^1_{\cV_1})   & H(K_1,\cT^1_{\cV_1\bcap f^*(\cV_2)}) \ar[l]_-{H(\kappa_1)} \ar[r]^-{H(\lambda_1)} &
  H(K_2,\cT^2_{\cV_2})   & \ldots \ar[l]_-{H(\kappa_2)}
}\\
\xymatrix{
  \ldots  & H(K_{n-1},\cT^{n-1}_{\cV_{n-1}\bcap f^*(\cV_n)}) \ar[l]_-{H(\kappa_n)} \ar[r]^-{H(\lambda_n)} &
   H(K_n,\cT^n_{\cV_n}) . 
}
\end{multline}
We refer to the persistence diagram of this module as the {\em persistence diagram of Morse decompositions}
of the sequence of combinatorial multivector fields\ $\cV_i$.


\section{Computational considerations and a geometric interpretation}
\label{sec:geom-int}
In this section we discuss computational aspects of the theory
and provide a geometric interpretation of the Alexandrov topology
of subsets of a simplicial complexes.

\subsection{Computational considerations}
Singular homology is not very amenable to computations.
Therefore, to compute the persistence module (possibly zigzag) in~\eqref{eq:persistence-module-simp} and~\eqref{eq:fi-zigzag-diagram1} efficiently, we take a more combinatorial approach.
We take the help of Theorem~\ref{thm:mccord}
(McCord's Theorem) in order to convert \eqref{eq:persistence-module-simp}
and \eqref{eq:fi-zigzag-diagram1} to a persistence module
where the objects are simplicial homology groups.

Let $(X,\cT)$ be a finite $T_0$ topological space and let $\leq_\cT$ be the partial order associated with $\cT$
by  Theorem~\ref{thm:alexandroff} (Alexandrov).
The {\em nerve} of this partial order,
that is, the collection of subsets linearly ordered by $\leq_\cT$ called {\em chains}, forms an abstract simplicial complex.
We denote it $N(X,\cT)$ or briefly $N(X)$ if $\cT$ is clear from the context.
Also by Alexandrov Theorem, a continuous map $f:(X,\cT)\to (X',\cT')$ of two finite topological $T_0$~spaces
preserves the partial orders
$\leq_\cT$ and $\leq_\cT'$. Therefore, it induces a simplicial map
$N(f): N(X,\cT)\to N(X',\cT')$.
Recall that every continuous and hence simplicial map $f: K\to K'$ of simplicial complexes extends linearly to a
continuous map $|f|: |K|\to |K'|$ on the polytopes of $K$ and $K'$ (cf. \cite[Lemma 2.7]{Munkres1984}).
The following proposition is straightforward.

\begin{prop}
\label{prop:barycentric-subdivision}
If $K$ is a simplicial complex, then the barycentric subdivision (cf. \cite[Sec. 2.15]{Munkres1984})
of a geometric realization of $K$ is a geometric realization of $N(K)$. In particular, $|K|=|N(K)|$.
Moreover, if $f:K\to K'$ is continuous, then $|f|=|N(f)|$.
\end{prop}

Consider the map
$
\mu_{(X,\cT)}:|N(X,\cT)|\ni x \mapsto \min\sigma_x\in X,
$
where $\sigma_x$ denotes the unique simplex $\sigma\in N(X,\cT)$ such that $x\in\cell{\sigma}$
and the minimum is taken with respect to the partial order $\leq_\cT$.
\begin{thm}
  \label{thm:mccord}
  (M.\ C.\ McCord, \cite{McCord1966})
  The map $\mu_{(X,\cT)}$ is continuous and a weak homotopy equivalence.
  Moreover, if $f:(X,\cT)\to (X',\cT')$ is a continuous map of two finite $T_0$ topological spaces, then
  the following diagrams commute.
  \begin{eqnarray*}
    \label{eq:mu-f-diag}
    \xymatrix{
      |N(X,\cT)| \ar[d]^{\mu_{(X,\cT)}} \ar[r]^{|N(f)|} & |N(X',\cT')| \ar[d]^{\mu_{(X',\cT')}}\\
      (X,\cT) \ar[r]^{f} & (X',\cT')} &
    \xymatrix{
      H_k(|N(X,\cT)|) \ar[d]^{H(\mu_{(X,\cT)})} \ar[r]^{H(|N(f)|)}
      & H_k(N(|X',\cT'|)) \ar[d]^{H(\mu_{(X',\cT')})}\\
      H_k(X,\cT) \ar[r]^{H(f)} & H_k(X',\cT')}
  \end{eqnarray*}
\end{thm}

By McCord's Theorem above,
there is a continuous map
$
\mu_{(X,\cT)}:|N(X,\cT)|\to X
$
which induces an isomorphism $H(\mu_{(X,\cT)}):H(|N(X,\cT)|)\to H(X,\cT)$ of singular homologies.
Moreover, the map $(X,\cT)\mapsto H(\mu_{(X,\cT)})$ is a natural transformation, that is
for any  continuous map $f:(X,\cT)\to (X',\cT')$  of finite $T_0$ topological spaces
$
H(\mu_{(X',\cT')})\circ H(|N(f)|) = H(f)\circ H(\mu_{(X,\cT)}).
$
Applying McCord's Theorem to every homology group in \eqref{eq:persistence-module-simp} we obtain the following
proposition.
\begin{prop}
\label{prop:nerve-morse-persistence}
Persistence module \eqref{eq:persistence-module-simp} is isomorphic to the persistence module
\begin{equation}
\label{eq:nerve-morse-persistence}
\xymatrix{
  H(|N(\union{\cM_1},\cT^1_{\cM_1})|) \ar[r]^-{f_1^N} &
  H(|N(\union{\cM_2},\cT^2_{\cM_2})|)\ar[r]^-{f_2^N} &
  \dots
  ~~H(|N(\union{\cM_n},\cT^n_{\cM_n})|),
}
\end{equation}
where $f_i^N:=H(|N(\bar{f}_i)|)$.
\end{prop}
Persistence module \eqref{eq:nerve-morse-persistence} is not yet simplicial,
but the map which sends each simplex in $K$ to the associated linear singular simplex in $|K|$
induces an isomorphism between the simplicial homology of $K$ and singular homology of $|K|$.
Moreover, this isomorphism commutes with the maps induced in simplicial and singular homology by simplicial maps
(see \cite[Theorems 34.3, 34.4]{Munkres1984}).
Thus, we obtain the following corollary. It facilitates the algorithmic computations of persistence diagrams for Morse decompositions of combinatorial dynamical systems.
\begin{cor}
\label{cor:simplicial-morse-persistence}
The persistence diagram of \eqref{eq:persistence-module-simp} is the same as the persistence diagram of the persistence module
\begin{equation}
\label{eq:simplicial-morse-persistence}
\xymatrix{
  H^\triangle(N(\union{\cM_1},\cT^1_{\cM_1})) \ar[r]^-{f_1^\triangle} &
  H^\triangle(N(\union{\cM_2},\cT^2_{\cM_2}))\ar[r]^-{f_2^\triangle} &
  \dots
  ~~H^\triangle(N(\union{\cM_n},\cT^n_{\cM_n})),
}
\end{equation}
where $H^\triangle$ denotes simplicial homology and $f_i^\triangle:=H^\triangle(N(\bar{f}_i))$. Moreover, an analogous statement holds for the zig-zag persistence module
\eqref{eq:fi-zigzag-diagram1}.
\end{cor}

For computing the persistence diagram of the module in ~\eqref{eq:simplicial-morse-persistence},
we identify the Morse sets in linear time by computing strongly connected components in $G_{F_i}$.
The nerve of these Morse sets can also be easily computed in time linear in input mesh size
(assuming the dimension of the complex to be constant).
Finally, one can use the persistence algorithm in~\cite{DFW14}, specifically designed for computing the persistence diagram of simplicial maps
that take the simplices of the nerve to the adjacent complexes in the sequence \eqref{eq:nerve-morse-persistence}.

\subsection{Geometric interpretation}
Proposition~\ref{prop:nerve-morse-persistence} provides means to interpret the Alexandrov topology of subsets of simplicial complexes
in the persistence module of Morse decompositions by the metric topology of their solids in the Euclidean space.
Recall that the solid of a subset $A\subset K$ of a simplicial complex is $|A|:=\bigcup\setof{\cell{\sigma}\mid\sigma\in A}$.
Let simplicial complexes $K_i$, combinatorial dynamical systems $F_i$ and Morse decompositions $\cM_i$ for $i=1,2,\ldots n$ be such as in Section~\ref{sec:IIS-MorseDecomp}. Moreover, assume $f_i:K_i\to K_{i+1}$ for $i=1,2,\ldots n$ are simplicial maps.
Let $\cO^i$ denote the metric topology of the polytope $|K_i|$.
Denote by $\famofgr{\cM}_i:=\setof{|M|\mid M\in\cM_i}$ the family of solids of Morse sets in $\cM_i$.
Consider the map $\nu_i:\union{\famofgr{\cM}_i}\ni x\mapsto |f_i|(x)\in \union{\famofgr{\cM}_{i+1}}$,
which is continuous with respect to topologies $\cO^i_{\famofgr{\cM}_i}$ and $\cO^{i+1}_{\famofgr{\cM}_{i+1}}$.

\begin{thm}
\label{thm:geometric-interpretation}
The persistence diagram of \eqref{eq:persistence-module-simp} is the same as the persistence diagram of the persistence module
\begin{equation}
\label{eq:geometric-interpretation}
\xymatrix{
  H(\union{\famofgr{\cM}_1},\cO^1_{\famofgr{\cM}_1}) \ar[r]^-{H( \nu_1)} &
  H(\union{\famofgr{\cM}_2},\cO^2_{\famofgr{\cM}_2})\ar[r]^-{H( \nu_2)} &
  \dots \ar[r]^-{H( \nu_{n-1})} & H(\union{\famofgr{\cM}_n},\cO^n_{\famofgr{\cM}_n}).
}
\end{equation}
\end{thm}
\proof
  By Proposition~\ref{prop:nerve-morse-persistence} it suffices to prove that the diagrams of \eqref{eq:nerve-morse-persistence} and \eqref{eq:geometric-interpretation} are isomorphic. By Theorem~\ref{thm:conn-comp-intersec}(iii), any two Morse sets
  in $\cM_i$ are disconnected.
  Hence, it follows from Proposition~\ref{prop:connected-finite}
  that the nerve $N(\union{\cM_i},\cT^i_{\cM_i})$
  splits as the disjoint union $\bigcup_{M\in\cM_i}N(M,\cT^i_M)$. In consequence, the whole diagram \eqref{eq:simplicial-morse-persistence}
  splits as the direct sum of diagrams for individual Morse sets.
  Again by Theorem~\ref{thm:conn-comp-intersec}(iii), any two sets in $\famofgr{\cM}_i$
  are $\cO^i_{\famofgr{\cM}_i}$-disconnected.
  Therefore, diagram \eqref{eq:geometric-interpretation} also splits as the direct sum of diagrams for individual sets  in $\famofgr{\cM}_i$.
  Thus, it suffices to prove that the respective diagrams for individual Morse sets are isomorphic.
  This follows easily from Proposition~\ref{lem:homotopy-equivalence} below. $\Box$\\

Note that, without loss of generality, given a simplicial complex $K$,  we may fix a geometric realization of $K$ and take its barycentric
subdivision as the geometric realization of  $N(K)$. Then, for any set of simplices $M\subset K$
we have $|N(M)|\subset |M|$.

\begin{prop}
\label{lem:homotopy-equivalence}
  The inclusion map $\iota_M:|N(M)|\to |M|$ is a homotopy equivalence.
  Moreover, if $f:K\to K'$ is a simplicial map, then the  map $\iota_M$
  and the map $\iota_{f(M)}:|N(f(M))|\to |f(M)|$ commute with the restrictions $|N(f)|_{\mid|f(M)|}$ and $|f|_{\mid |M|}$, that is
  $\iota_{f(M)}\circ |N(f)|_{\mid|f(M)|} = |f|_{\mid |M|}\circ \iota_M$.
  \label{lem:deform-retract}
\end{prop}

\proof To prove that $\iota_M$ is a homotopy equivalence, it suffices to show that $|N(M)|$ is
a deformation retract of $|M|$. To this end, order the simplices $\sigma_1,\ldots, \sigma_n$ in
$\cl M\setminus M$
 so that if $\sigma_j\preceq \sigma_k$, then $k \leq j$. Let $M_i=\cl M\setminus \{\sigma_1,\ldots,\sigma_i\}$ and consider the sequence
$
\cl M=M_0 \supseteq M_1 \supseteq\ldots \supseteq M_n=M.
$
We prove by induction on $n$ that $|M|$ deformation retracts to $|N(M)|$. Observe that the poset nerve $N(M_0)=N(\cl M)$
coincides with the barycentric subdivision of $\cl M$ and thus $|M_0|=|\cl M|=|N(M_0)|$.
Therefore, for $n=0$, the claim is satisfied trivially.

Inductively assume that $|M_{i-1}|$ deformation retracts to $|N(M_{i-1})|$ for all $i \in [1,n]$. We observe the following:

(1): In general $N(M_i)=N(M_{i-1})\setminus  C(\sigma_i)$ where $C(\sigma_i)$ denotes the set of all chains containing
$\sigma_i$ in the poset $(N(M_{i-1}),\preceq)$. If $\sigma_i^*$ denotes the vertex corresponding to $\sigma_i$ in $N(M_{i-1})$,
then $C(\sigma_i)$ is the star $\St\sigma_i^*$ in $N(M_{i-1})$. Also, $|\sigma^*|$ is the barycenter $b(\cell{\sigma})$.

(2): Let $Y \subset \St\sigma_i^*$ be any set of simplices in $N(M_{i-1})$ including $\sigma_i^*$. Then, $|N(M_{i-1})|\setminus |Y|$ deformation retracts to $|N(M_i)|$.
This follows from the fact that $|\St\sigma_i^*| \setminus |\sigma_i^*|=|\St\sigma_i^*| \setminus b(\cell{\sigma_i})$ retracts to the link of
 $\sigma_i^*$  along the segments that connect $\sigma_i^*$  to the points in the link and the restriction of this retraction to points in $\St\sigma_i^* \setminus |Y|$ provides the necessary deformation retraction.

For induction, observe that $|N(M_{i-1})|$ contains a subdivision of $|\sigma_i|=\cell{\sigma_i}$ because $M_{i-1}$ contains $\sigma_i$ and all its faces by definition of $M_i$s.
Let $Y$ denote the set of simplices that subdivide $\cell{\sigma_i}$. Then, according to (2), $|N(M_{i-1})|\setminus \cell{\sigma_i}$ deformation retracts to $|N(M_i)|$. We construct a deformation retraction of $|M_i|$ to $|N(M_i)|$ by first retracting  $|M_{i-1}|$ to $|N(M_{i-1})|$ by the inductive hypothesis and then retracting $|N(M_{i-1})|\setminus \cell{\sigma_i}$ to $|N(M_i)|$.
The remaining part of the lemma is an immediate consequence of Proposition~\ref{prop:barycentric-subdivision}. $\Box$


\begin{figure*}[ht!]
  \setlength\tabcolsep{1.5pt}
  \centering
  \includegraphics[width=0.33\textwidth]{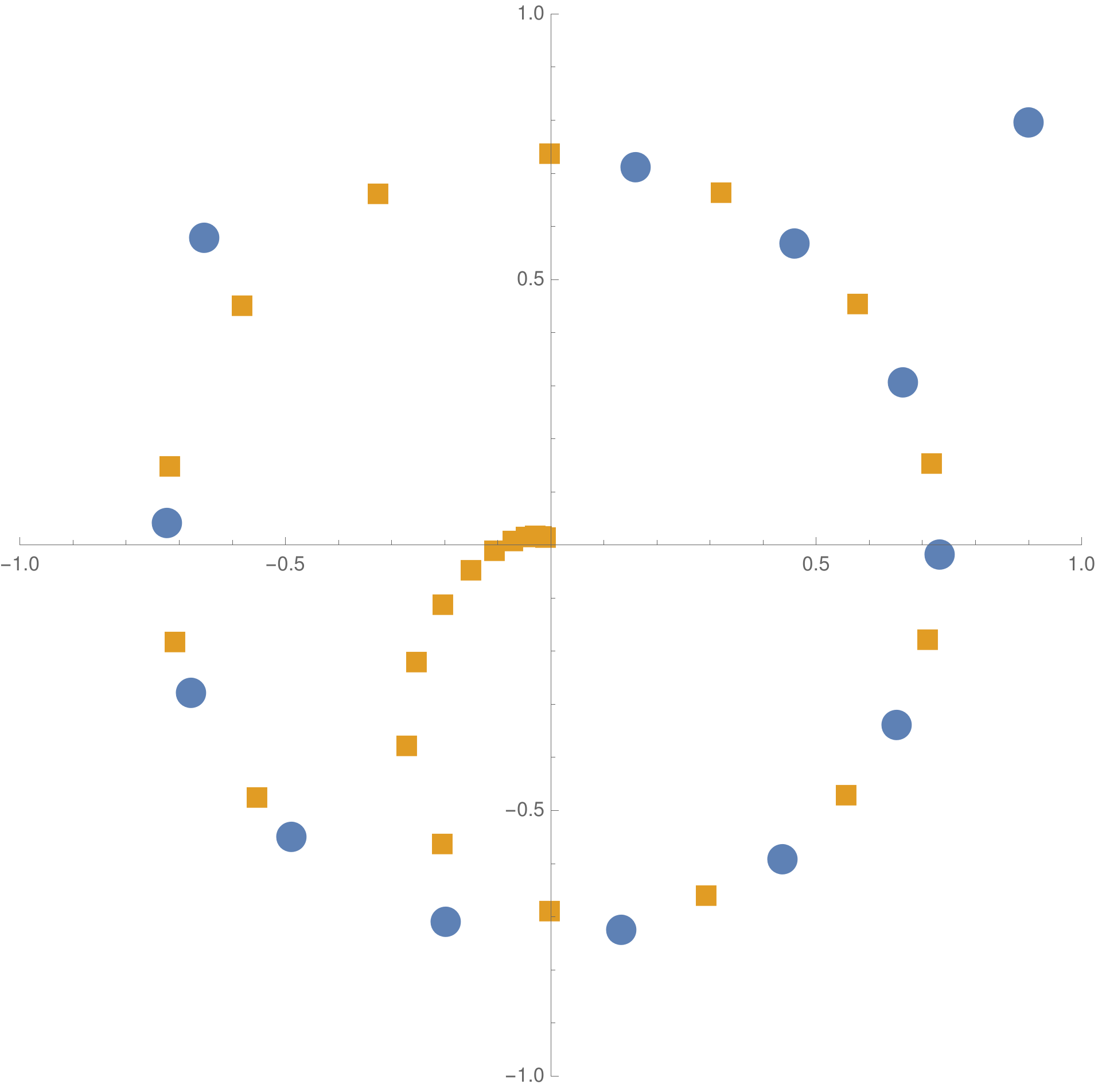}
  \includegraphics[width=0.33\textwidth]{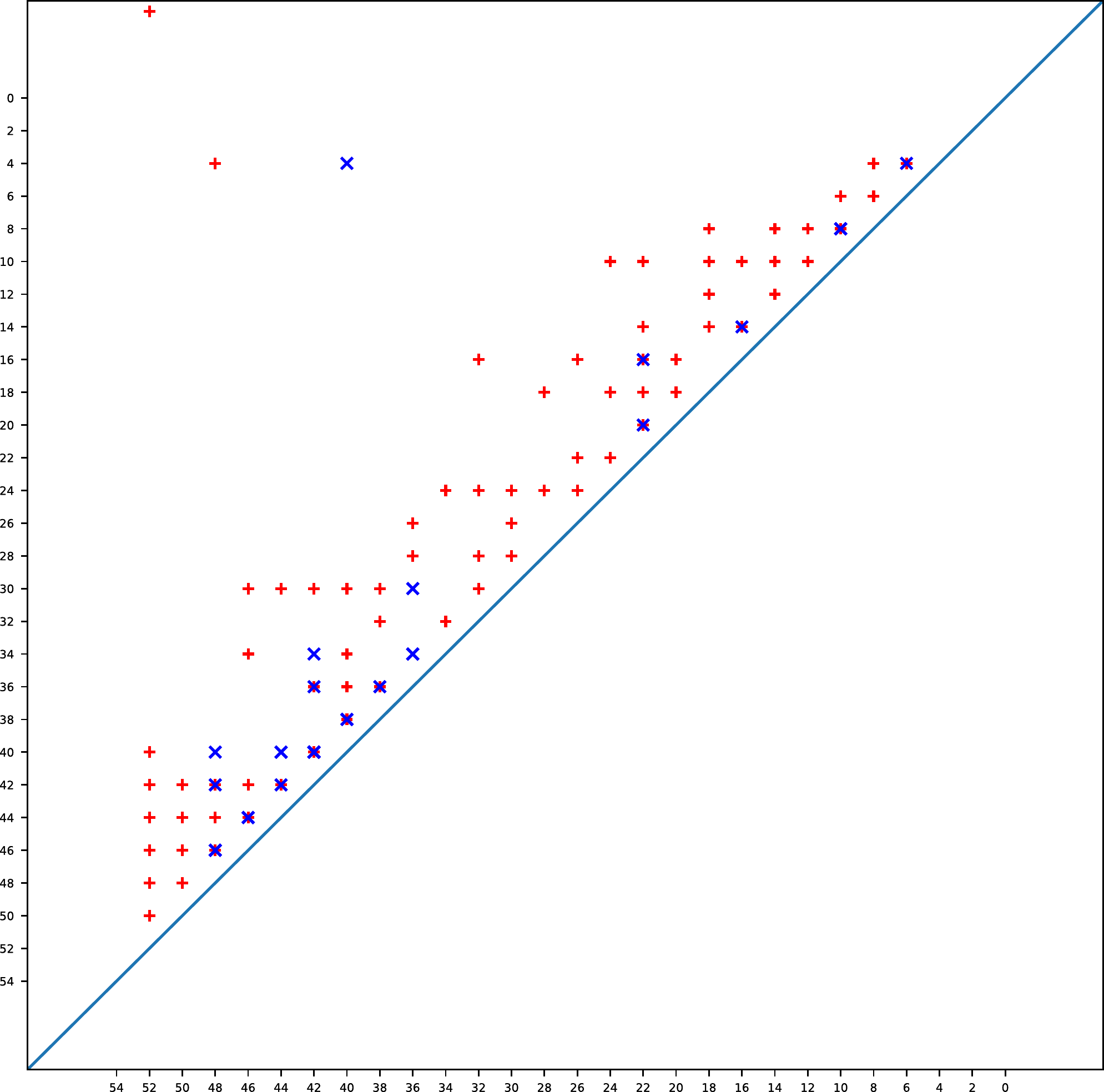}\\
  \medskip
  \includegraphics[width=0.24\textwidth]{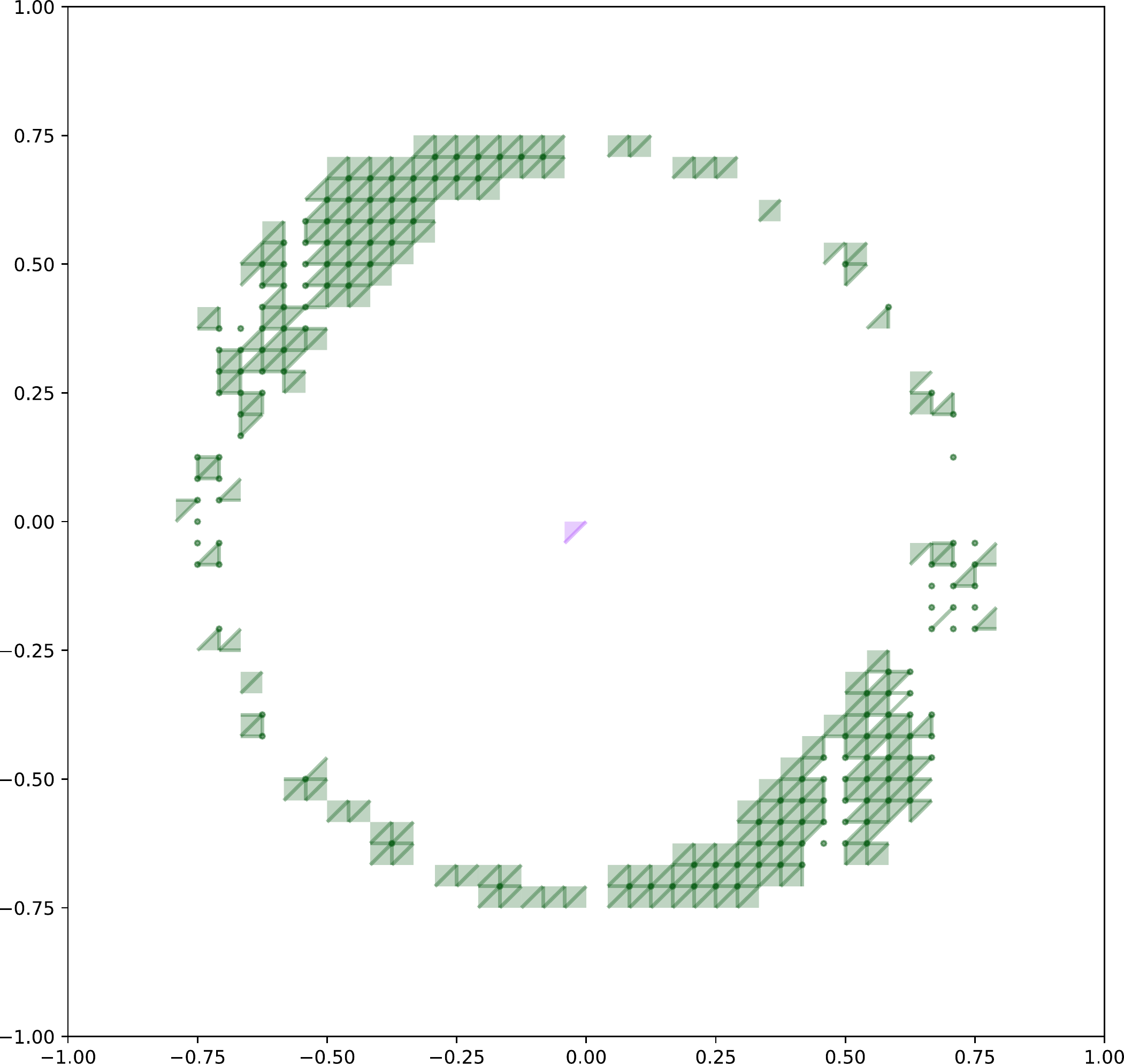}
  \includegraphics[width=0.24\textwidth]{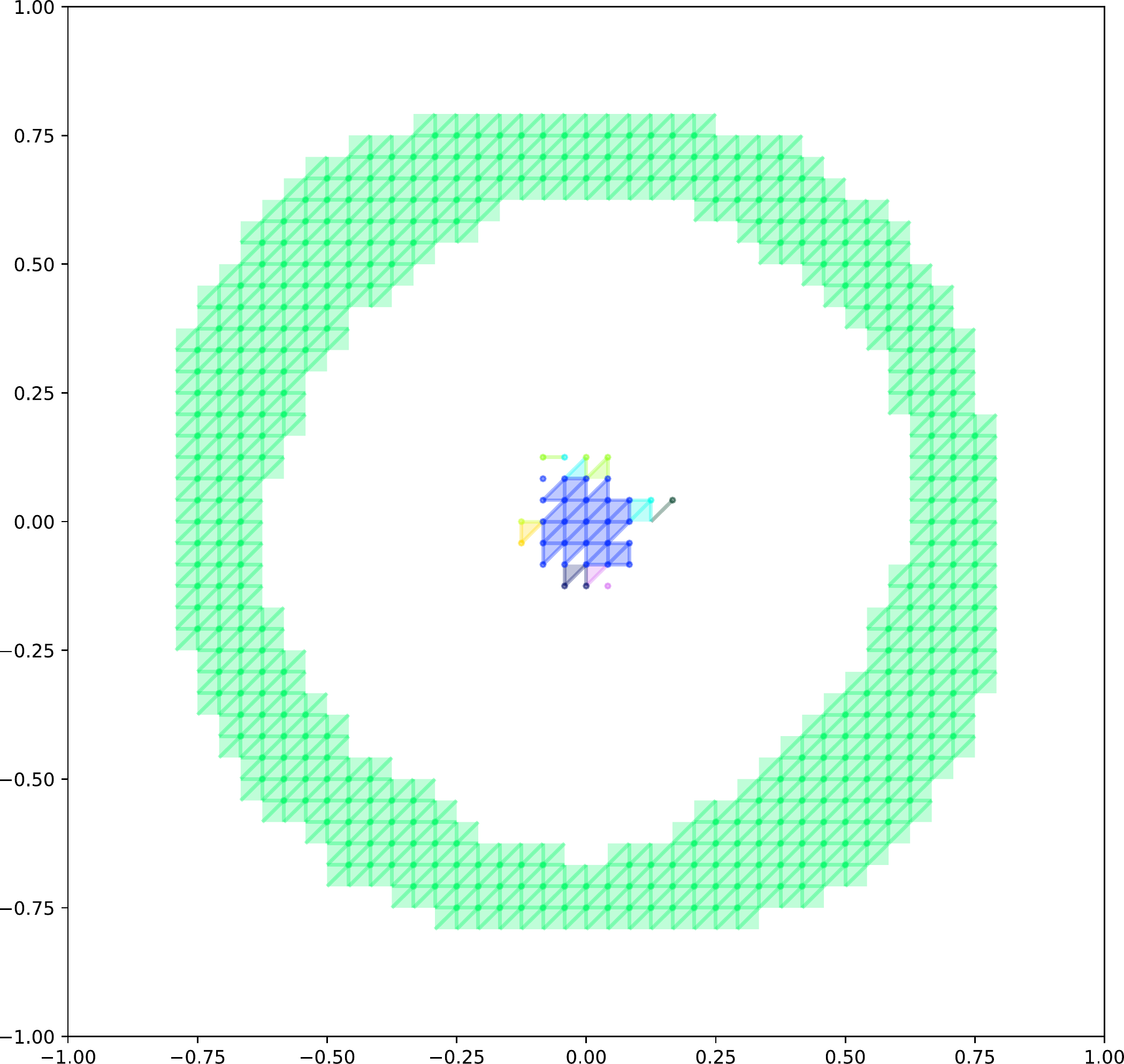}
  \includegraphics[width=0.24\textwidth]{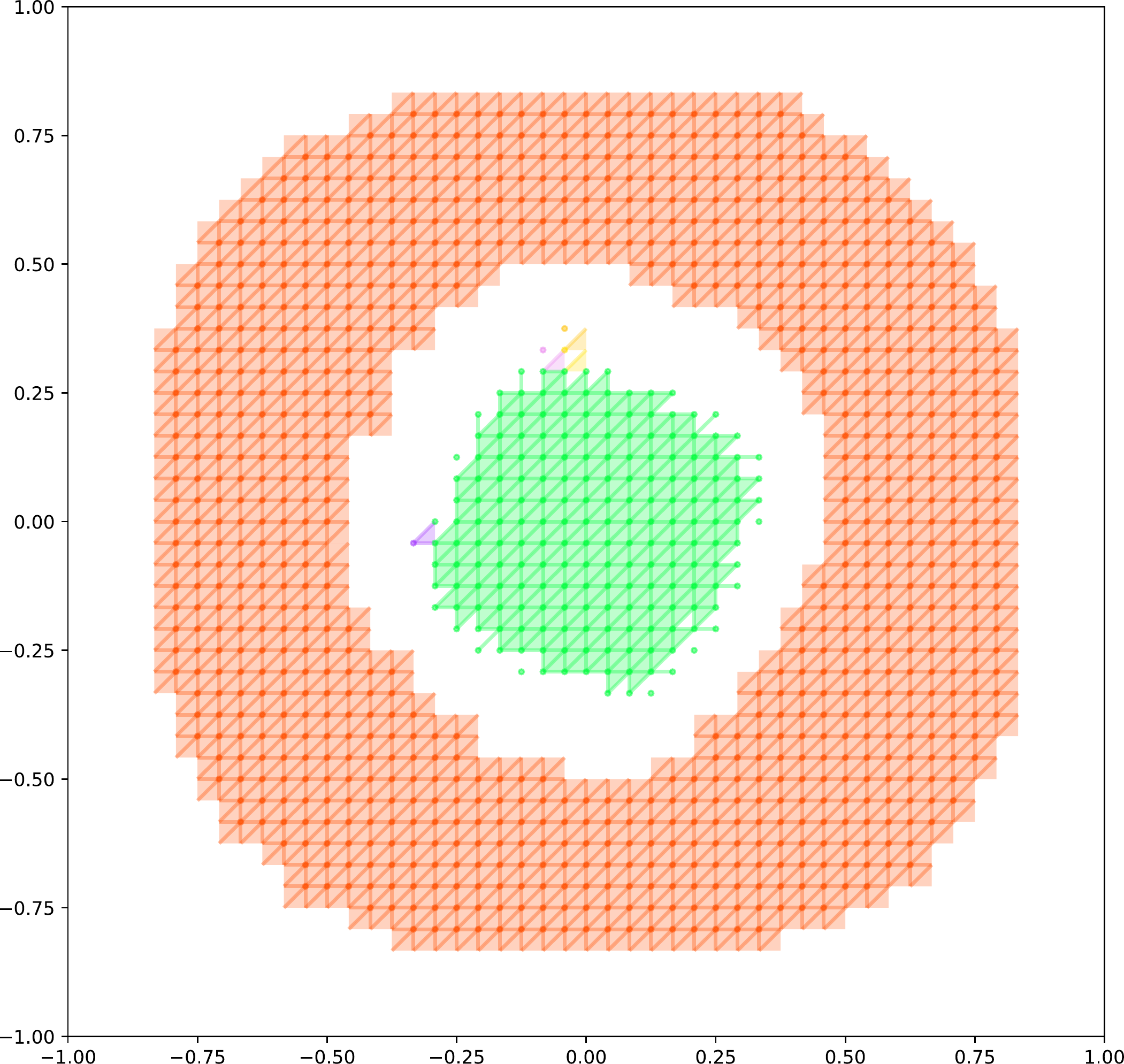}
  \includegraphics[width=0.24\textwidth]{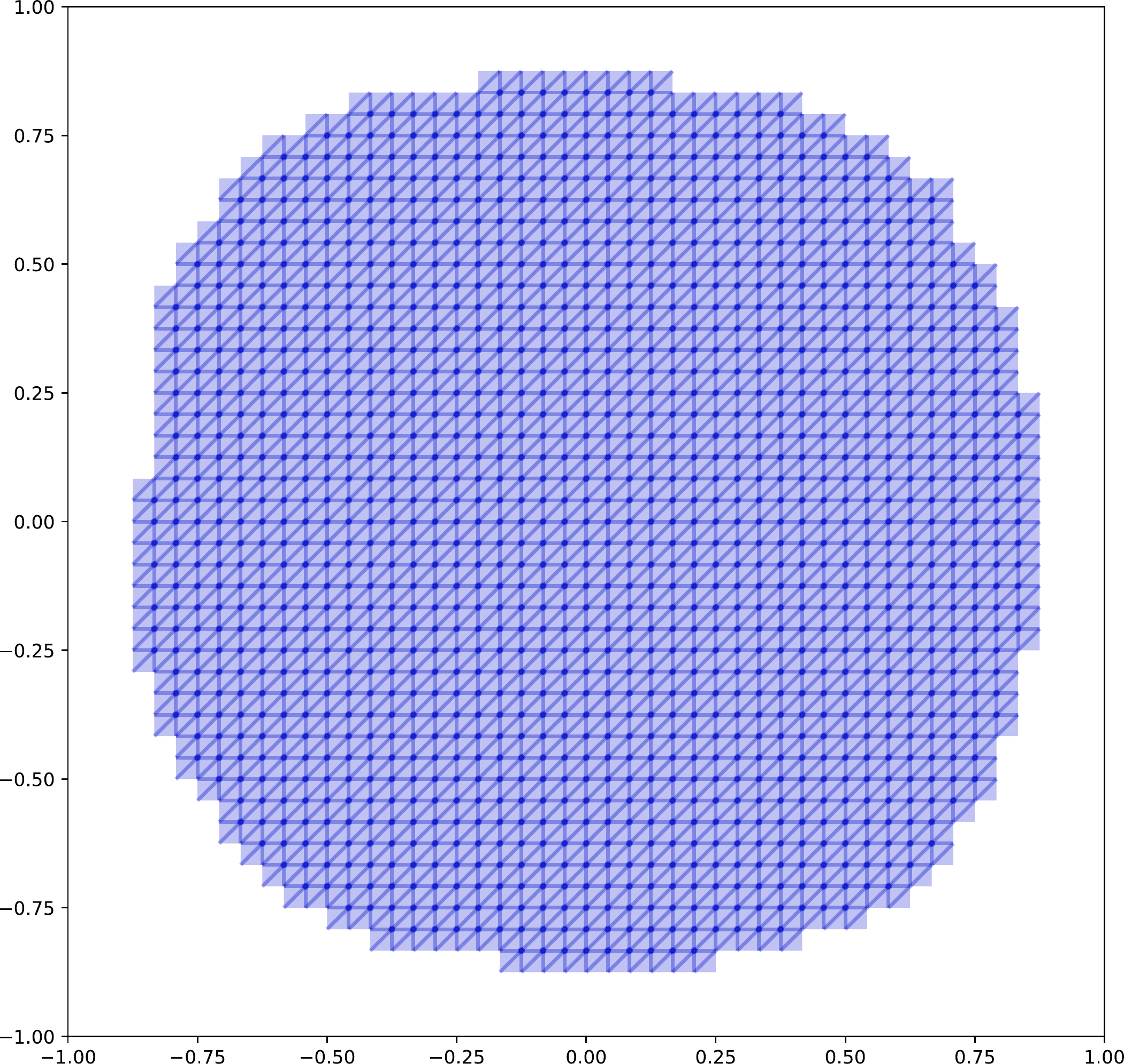} \\
  \caption{
    Upper left: two trajectories of  (\ref{eq:neimark_sacker_map}) with starting points at $(-0.01,0.01)$ (squares) and $(0.9, 0.8)$ (circles).
    Upper right: persistence diagram for a noisy sampling of \eqref{eq:neimark_sacker_map}.
    Red pluses and blue crosses indicate homology generators in dimension zero and one, respectively.
    Bottom from left to right: Morse sets in selected filtration steps for threshold values
    $\mu=\frac{48}{148}$ 
    , $\mu=\frac{30}{148}$ 
    , $\mu=\frac{6}{148}$ 
    and $\mu=\frac{4}{148}$ 
    , respectively.
  }
  \label{fig:ex_2d_map_persistence}
\end{figure*}

\section{Examples}
\label{sec:examples}
In this section we present two numerical examples.
The first example concerns the persistence of
the Morse decompositions of a noisy sample of
Kuznetsov map with respect to a frequency parameter.
The second example concerns the persistence of
the Morse decompositions of combinatorial multivector fields
with respect to an angle parameter of the algorithm constructing
the fields from a cloud of vectors.

\subsection{Kuznetsov map.}
\label{subsec:kuzn_ex}
Let us consider the following planar map analyzed by Kuznetsov in the context of
the Neimark-Sacker bifurcation ~\cite[Subsection 4.6]{Kuzn95}.
\begin{equation}
\label{eq:neimark_sacker_map}
N\bigg(\begin{bmatrix} x_{1} \\ x_{2} \end{bmatrix}\bigg) :=
    \begin{bmatrix} \cos(\theta) & -\sin(\theta) \\ \sin(\theta) & \cos(\theta) \end{bmatrix}
     \bigg((1+\alpha)\begin{bmatrix} x_{1} \\ x_{2} \end{bmatrix} + (x_{1}^2 + x_{2}^2)
     \begin{bmatrix} a & -b \\ b & a \end{bmatrix} \begin{bmatrix} x_{1} \\ x_{2} \end{bmatrix}\bigg).
\end{equation}

For parameters $\theta = \pi/17$, $\alpha=0.5$, $a=-1$ and $b=0.5$
the system restricted to square $[-1, 1]\times[-1, 1] \subset \mathbb{R}^2$
admits a Morse decomposition consisting of an unstable fixed point and an attracting invariant circle.
(see Figure \ref{fig:ex_2d_map_persistence}, upper left).
We want to detect this Morse decomposition
just from a finite sample of the map and in the presence of  Gaussian noise.
The setup is similar to the toy example in Section \ref{sec:comb-dyn}.

Let $x\in\RR^2$ and $\epsilon_X,\epsilon_Y\in\RR^2$ be random vectors
chosen from normal distribution centered at zero,
with standard deviation $\sigma_X$ and $\sigma_Y$ respectively.
Let
\begin{equation}
\tilde{N}(x) := N(x+\epsilon_X)+\epsilon_Y
\end{equation}
be a noisy version of the map (\ref{eq:neimark_sacker_map}).
Consider a triangulation $K$ of the square  $Q:=[-1, 1]\times[-1, 1]\subset \mathbb{R}^2$
obtained by splitting $Q$ into a $48\times 48$ uniform grid of squares of size $r=1/24$ and dividing every square into two triangles.
Then, the set of toplexes $K_{top}$ consists of $4608$ 2-simplices.
The $\sigma_X$, $\sigma_Y$ are taken to be proportional to the grid size $r$, that is $\sigma_X=r/4$ and $\sigma_Y=r$.
A noisy sample $\Gamma=\{(x_i,y_i)\}_i$ of the map $N$ is generated by taking an uniformly distributed sequence of points $x_i$ in $Q$
and its disturbed images $y_i:=\tilde{N}(x_i)$.
Pairs $(x_i,y_i)\in\Gamma$ such that $y_i\not\in Q$ has been rejected from a sample.
The combinatorial dynamical system $F_\mu$ is constructed in the same way as in Section \ref{sec:comb-dyn}, namely
\begin{equation}\label{eq:ex_2d_cds}
F_\mu(\sigma):=
 \co\bigcup_{\tau\in K_{top},\,\sigma\preceq\tau}\setof{\bar{\tau}\in K_{top}\mid
  \frac{n_{\tau,\bar{\tau}}}{n_{max}} \geq \mu },
\end{equation}
where $n_{\tau,\bar{\tau}}$ denotes the number of pairs in $\Gamma$ connecting two toplexes $\tau,\bar{\tau}\in K_{top}$, that is
\begin{align}
n_{\tau,\bar{\tau}}:=\#\{(x_i,y_i)\mid x_i\in \cl |\tau| \text{ and } y_i\in \cl |\bar{\tau}|\}
\end{align}
and $n_{max}$ is maximal of these values. For this particular experiment $n_{max}=148$. Note that construction of $F_\mu$  (\ref{eq:ex_2d_cds}) is also well defined for lower dimensional simplices. 

The parameter $\mu$ in (\ref{eq:ex_2d_cds}) describes the minimal frequency of an edge to be present in the combinatorial dynamical system $F_\mu$.
The family of Morse sets $\mathcal{M}(F_\mu)$ at given level $\mu$ consists of all strongly connected components of an associated graph. The set of considered frequency levels $1=\mu_{74} > \mu_{73} > ... > \mu_0=0$, where $\mu_i=\frac{2i}{148}$, leads to the sequence of Morse decompositions such that $\mathcal{M}(F_{\mu_i}) \sqsubseteq \mathcal{M}(F_{\mu_{i-1}})$. The persistence diagram at Figure \ref{fig:ex_2d_map_persistence}(upper right), for clarity, shows only results for $\mu\leq\mu_{27}=\frac{54}{148}$, since this is a level where Morse sets start to emerge.

Results are presented in Figure \ref{fig:ex_2d_map_persistence}.
As expected, the persistence diagram (Figure \ref{fig:ex_2d_map_persistence}, upper right) indicates the presence of two 0-dimensional and one 1-dimensional homology generators with high persistence.
Bottom row at Figure \ref{fig:ex_2d_map_persistence} shows Morse sets for selected frequency levels. For lower thresholds, both invariant sets eventually merge together creating the only strongly connected component. In the case without noise, fixed point at the origin and attracting invariant set remains separated for all values of $\mu$.

\begin{figure*}[ht!]
  \begin{center}
    \includegraphics[width=0.3\textwidth]{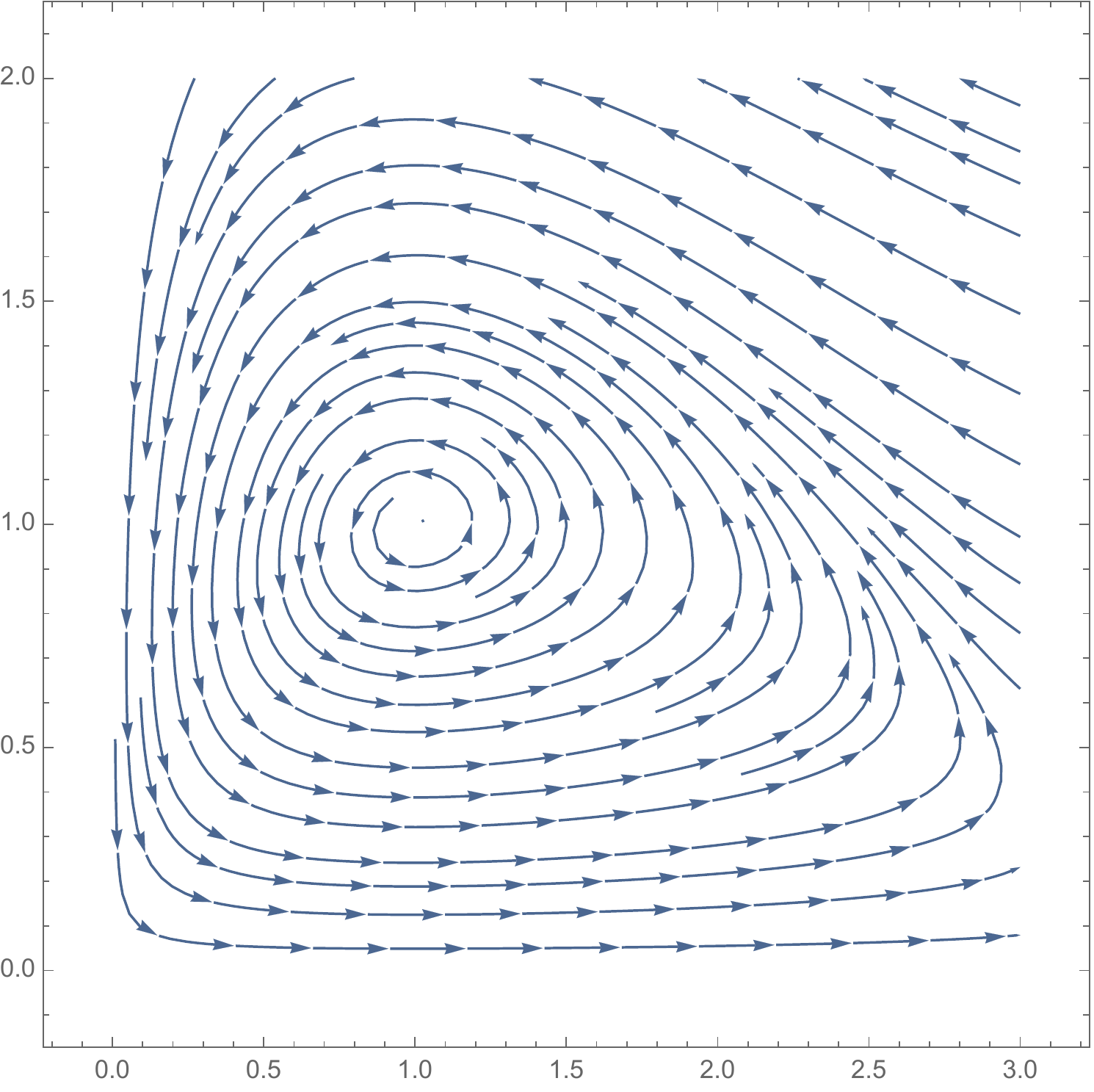}
    \includegraphics[width=0.3\textwidth]{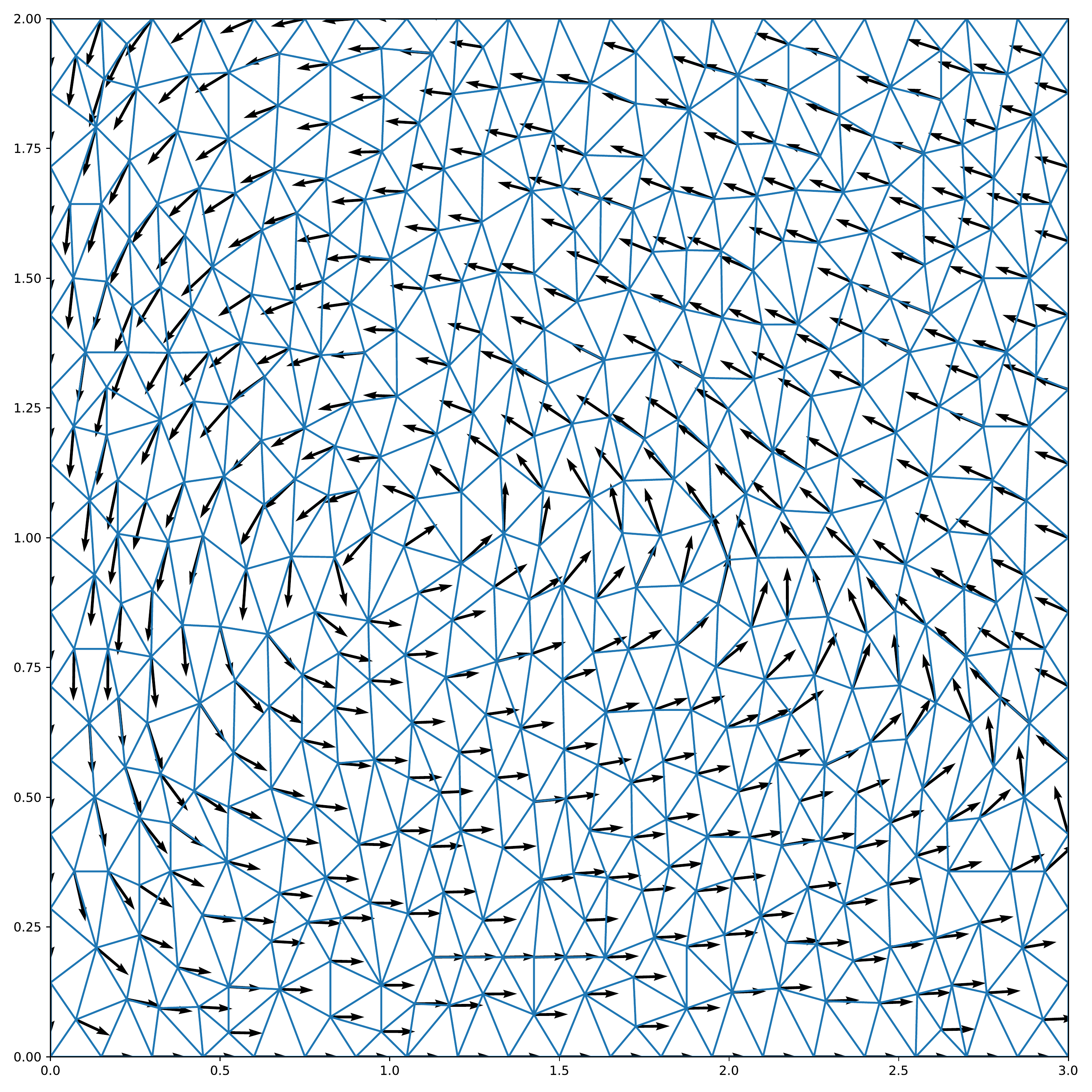}
  \end{center}
  \vspace{-0.7cm}
%
  \begin{center}
    \includegraphics[width=0.3\textwidth]{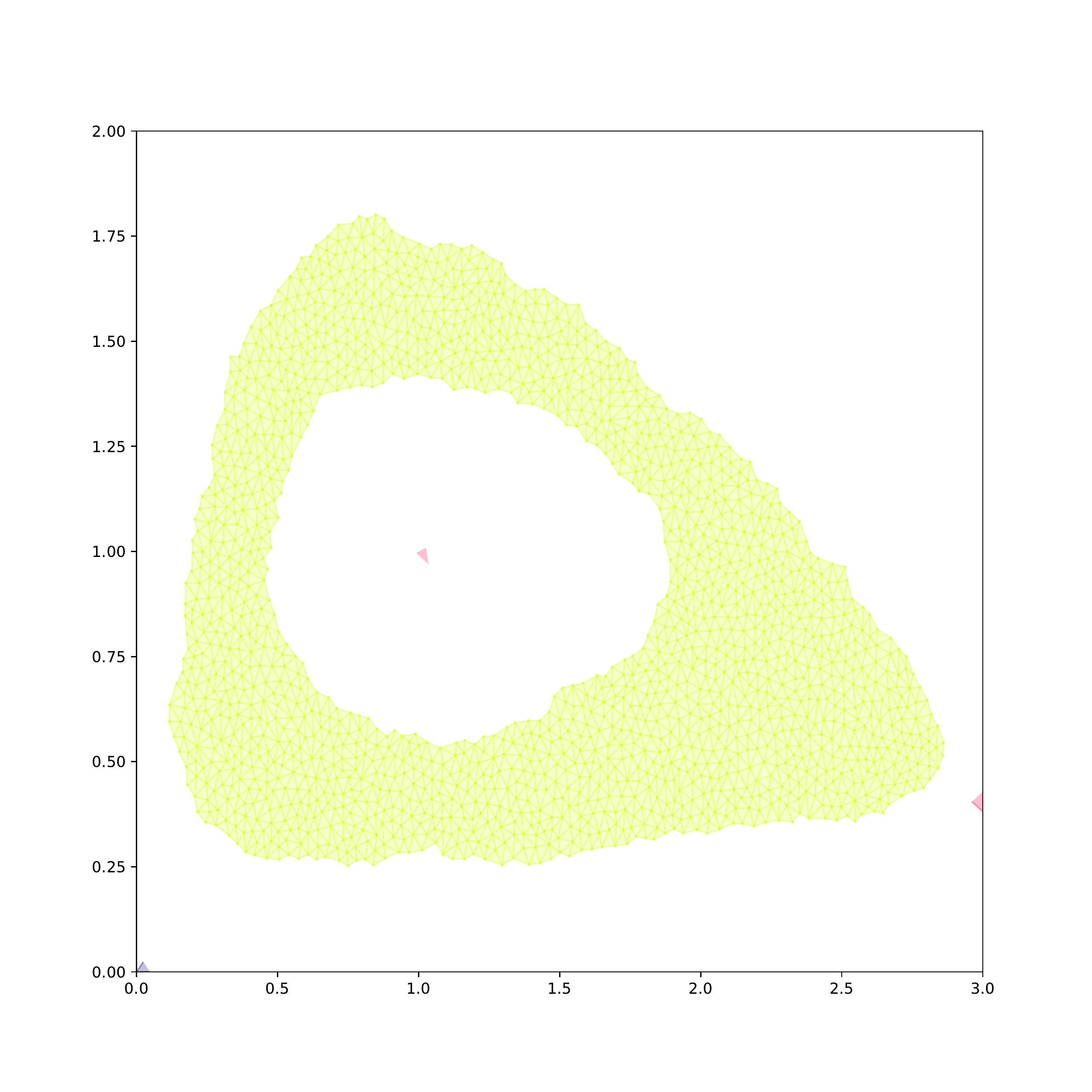}
    \includegraphics[width=0.3\textwidth]{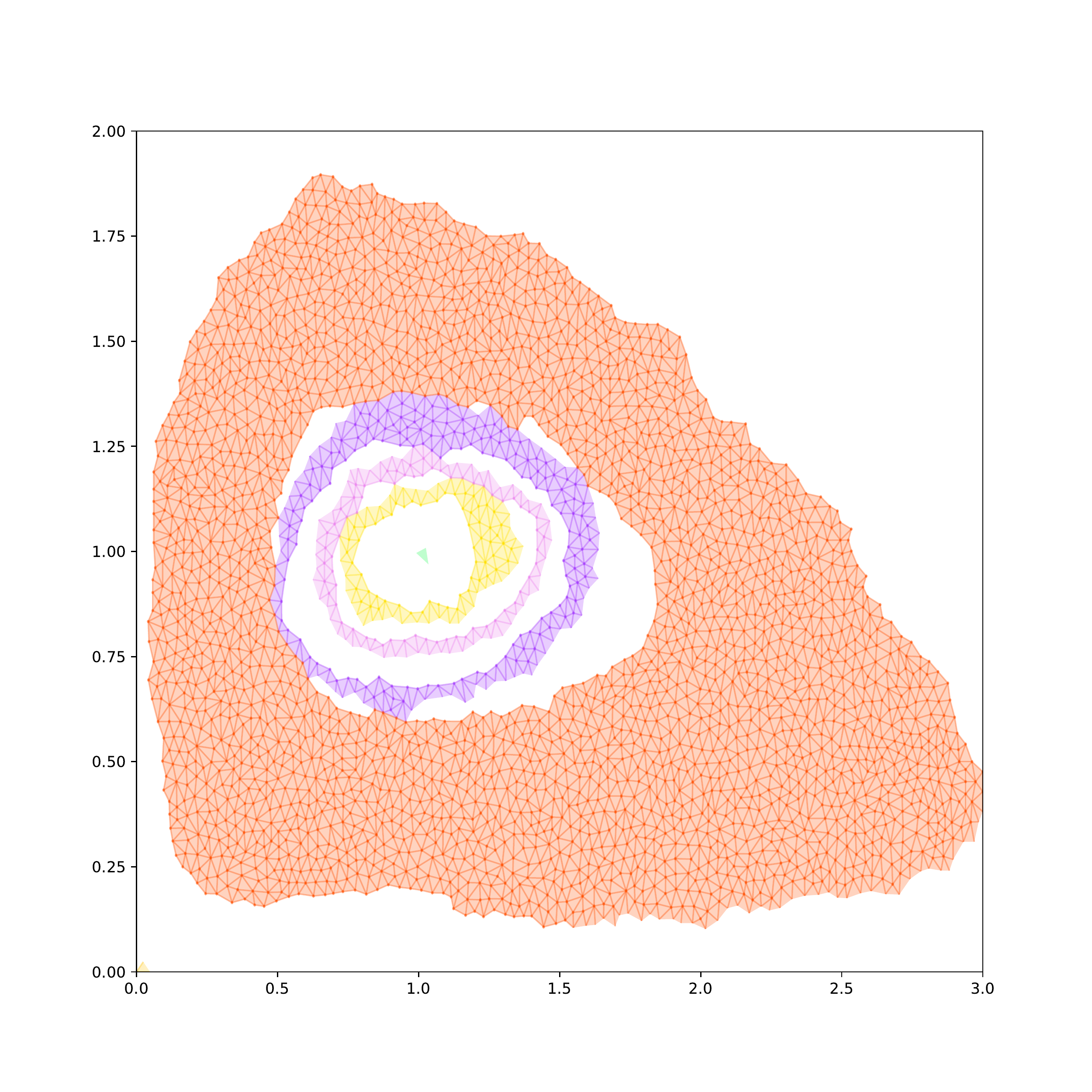}
    \includegraphics[width=0.3\textwidth]{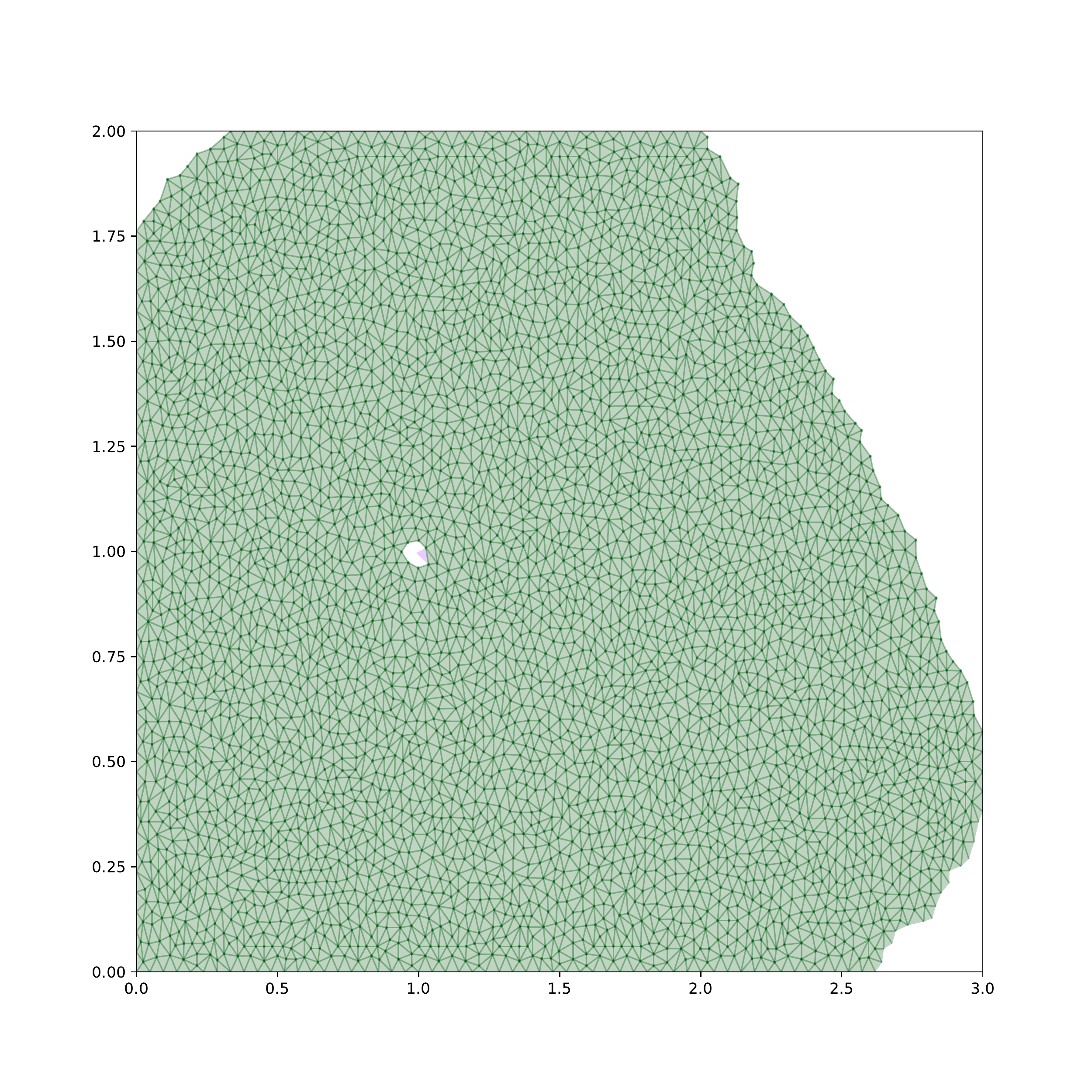}
  \end{center}
 \vspace{-1cm}
  \begin{center}
    \includegraphics[width=0.8\textwidth]{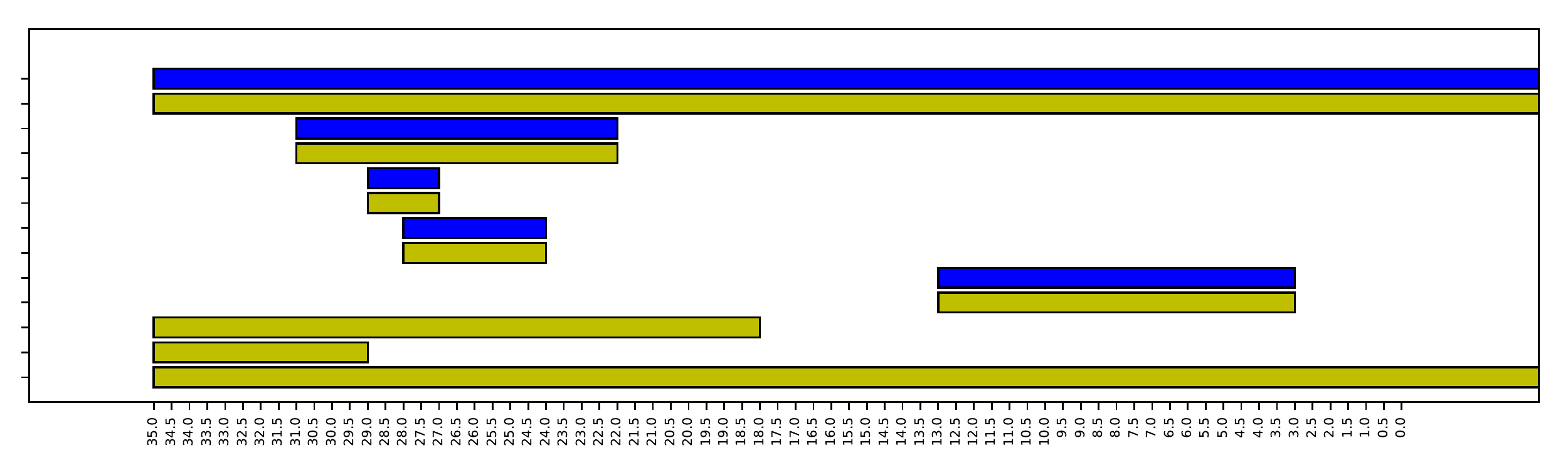}
  \end{center}
  \vspace{-0.7cm}
  \caption{
    Top: A streamline plot and a cloud of vectors on a triangular grid for LV model.
    Middle from left to right: Morse decompositions of combinatorial dynamical systems constructed from the cloud
    for $\alpha=35^{\circ}$, $\alpha=28^{\circ}$ and $\alpha=0^{\circ}$, respectively.
    The mesh used contains $4756$ points and $9300$ triangles.
    Bottom: The respective persistence barcode. Yellow bars represents generators in dimension $0$ and blue bars in dimension $1$.
    The bars are ordered by average size of the Morse sets (biggest on  the top).
  }
  \label{fig:lv-zigzag}
\end{figure*}
\subsection{Lotka Volterra model.}
\label{sec:exlv}
Consider the Lotka-Volterra (LV) model:
\begin{equation}
\frac{\partial x}{\partial t} =  x\left(1-\frac{x}{k}\right) - \frac{(a_1  x y  )}{b + x},  \qquad \frac{\partial y}{\partial t} =
\frac{a_2  x y}{b + x} - g  y,
\end{equation}
where $k = 3.5$, $b = 1$, $g = 0.5$, $a_1 = (1 - \frac{1}{k}) (b + 1)$, $a_2 = g  ( b + 1)$
(see ~\cite[Chapter 2, Eq. 2.13 and 2.14]{Boccara2004}).
The system has a Morse decomposition consisting of
a repelling stationary point
and an attracting periodic orbit.
We want to observe this Morse decomposition in a combinatorial dynamical system constructed from a finite sample of the vector field.
In Table~\ref{tab:algorithm}
we present an algorithm for constructing a combinatorial multivector field from a sampled vector field.
The algorithm requires an angle parameter $\alpha$. The constructed combinatorial multivector field and hence its
combinatorial dynamical system depend on this parameter.
We execute the algorithm for  varying $\alpha$
and construct the zigzag filtration \eqref{eq:fi-zigzag-diagram1}.
Since the supporting simplicial complex (mesh) remains fixed, we obtain zigzag persistence
under inclusion maps.
Experiments with varying mesh, utilizing non-inclusion maps, are in progress.
The outcome for the LV model is presented in Figure~\ref{fig:lv-zigzag}.
We note that the trivial Morse sets that is Morse sets consisting of just one multivector $V$
such that $H(\cl V,\cl V\setminus V)=0$ are excluded from the presentation of Morse decompositions and from the barcode, because
such Morse sets are considered spurious due to the triviality of their Conley index (see \cite{Mr2016}).

%
The input to the algorithm $\mathrm{CVCMF}$ in Table~\ref{tab:algorithm} that computes a multivector field from a cloud of vectors consists of:
\begin{itemize}
  \item[-] a simplicial mesh $K$ with vertices in a cloud of points $P=\{p_i\mid i=1,2,\ldots, n\}\subset \mathbb{R}^d$,
  \item[-] the associated cloud of vectors $V := \{   \vec{v_i} \mid i=1,2,\ldots, n\}\subset\RR^d$ such that vector $\vec{v_i}$ originates from point $p_i$,
  \item[-] an angular parameter $\alpha$.
\end{itemize}

\begin{table}[h!]
  \begin{algorithmic}[1]
    \Procedure{CVCMF}{$K$, $V$, $\alpha$}
    \State $m \gets$ an identity map $K \to K$.
    \ForAll{$\sigma \in K$}
    \State $m[\sigma] \gets$ any toplex in the star of $\sigma$ pointed by mean of $\setof{\vec{v_i} \in V \mid p_i \preceq \sigma}$
    \EndFor
    \ForAll{$ i=1,2,\ldots, n $} \Comment{Aligns vectors}
    \State $S \gets$ $\setof{(\dim \sigma, \measuredangle(\sigma,\vec{v_i}), \sigma) \mid \sigma \in K \text{ and } p_i \preceq \sigma \text{ and } \measuredangle(\sigma,\vec{v_i}) \le \alpha}$
    \State $S' \gets$ sort $S$ using lexicographical order on first two positions \Comment{$(\dim , \measuredangle , \_)$}
    \State $(\_,\_,\sigma) \gets$ first element of $S'$
    \State $m[p_i] \gets \sigma$
    \EndFor
    \ForAll{$\sigma \in K$ in descending dimension} \Comment{Remove convexity conflicts}
    \While{exists $\tau \preceq \sigma$ s.t. $[\tau, m[\tau]] \cap [\sigma, m[\sigma]] \ne \emptyset$ and $m[\tau] \neq m[\sigma]$}
    \State $m[\tau]\gets \sigma$ and $m[\sigma] \gets \sigma$
    \EndWhile
    \EndFor
    \State $\cV \gets$ build a partition of $K$ using nonempty pre-images of $m$
    \State \Return $\cV$.
    \EndProcedure
  \end{algorithmic}

  \vspace{0.1in}
  \caption{An algorithm constructing a combinatorial multivector field from a sampled vector field}
  \label{tab:algorithm}
\end{table}

For each simplex $\sigma\in K$ and a vector $\vec{v_i}$ originating from
vertex $p_i$ of $\sigma$
we measure the angle between $\vec{v_i}$ and the affine subspaces
spanned by the vertices of  $\sigma$.
We assume the angle to be zero when the vector has length zero
or the simplex is just a vertex.
For a toplex $\sigma$, we assume that the angle is zero when $\vec{v_i}$
points inward $\sigma$ and $\infty$ otherwise.
When the angle is smaller than $\alpha$, we project $\vec{v_i}$ onto $\sigma$.
Intuitively, it aligns the vectors to the lower dimensional simplices.
After this alignment, a multivector field is constructed by removing the convexity conflicts.
Obviously, the output depends on the parameter $\alpha$.
We measure changes in the multivector field $\cV$ via persistence of its Morse decomposition.
To compute such persistence we use Dionysus software~\cite{MorozovZigzag}.

\end{document}